%% file: attitude.tex
\documentclass[journal, letterpaper]{IEEEtran}
\usepackage[normalem]{ulem}

\ifCLASSINFOpdf
  \usepackage[pdftex]{graphicx}
  \graphicspath{{.}}
  \DeclareGraphicsExtensions{.pdf,.jpeg,.png}
\else
\fi

\usepackage[cmex10,reqno]{amsmath}
\newcommand{\vect}[1]{\mathbf{#1}}
\usepackage{xfrac}

\usepackage{tabularx}
\usepackage{multirow}
\usepackage{chngcntr}
\counterwithout{table}{section}

\usepackage{url}
\begin{document}
%
\title{Comparison of Attitude Estimation Techniques for Low-cost Unmanned Aerial Vehicles}
%
%
%

\author{Harris Teague, \textit{Qualcomm Research, Qualcomm Technologies, Inc.,} San Diego, California}
\maketitle

\begin{abstract}
Attitude estimation for small, low-cost unmanned aerial vehicles is often achieved using a relatively simple complementary filter that combines onboard accelerometers, gyroscopes, and magnetometer sensing.  This paper explores the limits of performance of such attitude estimation, with a focus on performance in highly dynamic maneuvers.  The complementary filter is derived along with the extended Kalman filter and unscented Kalman filter to evaluate the potential performance gains when using a more sophisticated estimator.  Simulations are presented that compare performance across a range of test cases, many where ground truth was generated from manually controlled flights in a flight simulator.  Estimator scenarios that are generic across the different estimator types (such as the way sensor information is processed, and the use of dynamically changing gains) are compared across the test cases.  An appendix is included as a quick reference for the common attitude representations and their kinematic expressions.

\end{abstract}


%


\section{Introduction}
%
%
%
%
\IEEEPARstart{A}{ttitude} estimation is a fundamental component of all vehicle estimation and control systems and there is a rich body of  literature on this subject.  However, implementation of the various techniques can be time-consuming and intricate.  Often, engineers face the question of what algorithm and architecture to choose for their application, but there is limited information available to help them make this decision.  The question comes down to, "will implementation of algorithm \textit{X} provide performance improvements that justify increased algorithm complexity?"  This paper attempts to help answer this question for the case of an attitude estimator that employs consumer-grade inertial and magnetometer components.  It also sets a framework that should be useful for others who wish to answer the same question for other cases.

This paper starts with a brief survey of the literature.  Next, equations are developed for what we view as the three primary algorithms in use today - the complementary filter (CF), the extended Kalman filter (EKF), and the unscented Kalman filter (UKF).  These are developed for the generic case of vector measurements, and includes some practical considerations for the use of consumer-grade sensing.  A goal is to provide a concise development of these nonlinear attitude estimation techniques such that the reader can quickly and easily implement working algorithms.  

Next, simulation results are presented for a range of scenarios and assumptions and the performance of the estimators is compared and contrasted, and some conclusions drawn.  Since the most challenging case for the estimators is during times of highly dynamic maneuvers (that produce significant translation accelerations), the simulations focus on such cases.  To produce realistic motions, test case ground truth was produced by recording the outputs of a flight simulator while the desired motions were flown manually by a pilot.  Also, since there are several different attitude representations that can be used, the appendix shows the kinematic equations for the most common representations.   

\subsection{Literature}

There are many papers and texts on attitude estimation in the literature, so the focus here will be on the work that was most influential to this paper.  A comprehensive summary of different attitude estimation methods is given in \cite{crassidis2007survey} which includes a full spectrum of approaches including filter/smoothers, Wahba's problem techniques, and 2-step methods.  However, here we focus on quaternion-based estimators that employ attitude perturbation states in the filter, and one of the important early papers that uses this approach is \cite{toda1969spars}.  These techniques are extended in \cite{lefferts1982kalman} and \cite{miller1983new}.  Convergence and stability properties are studied in \cite{salcudean1991globally, mahony2005complementary,nonlinear_comp_filt_mahony_08} and \cite{jensen2011generalized}.  The multiplicative extended Kalman filter (MEKF) is derived and justified in \cite{markley2002fast,markley2003attitude,markley2004attitude}. An analytical comparison of CF and EKF is given in \cite{Higgins1975} (without comparison of any results). UKF techniques introduced in \cite{julier1997new} are explored for attitude estimation in \cite{crassidis2003unscented, laviola2003comparison, perea2007nonlinearity}, and \cite{perea2008new}.  The attitude representations that are used in these references (and in this paper) are surveyed and explained in \cite{shuster1993survey}. 

\subsection{Common algorithm elements}

Straightforward development of an attitude estimator starts with defining a state vector
that fully specifies the time-varying system. For example, the attitude quaternion $q$ and
the attitude rate or angular velocity $\vect{\omega}$ can specify the system. Further, it is
common to extend this state to include measurement error terms that have components that are
not independent over time, such as slowly drifting biases and misalignments. Focusing just
on gyro bias, we can construct an attitude estimator with the following time-varying states:

\begin{equation}
	x = \begin{bmatrix} q \\ \omega \\ b\end{bmatrix},
\begin{array}{l} \text{atttitude unit quaternion} \\ \text{angular velocity} \\ \text{gyro
bias}\end{array}
\end{equation}

However, we will not use this state, and will instead do two things:

\begin{enumerate}
\item Use a 3-vector representation of attitude \cite[Section VI]{toda1969spars} such as Euler vector  $\vect{a}_\phi$, Gibbs vector  $\vect{a}_g$, or
modified Rodriques parameters  $\vect{a}_p$(see Appendix \ref{app:att}). Since q is unit magnitude and
redundant, its expected value is not always a valid quaternion, so is not theoretically a
good choice for linearized optimal estimation \cite{markley2003attitude}. As we will see, it
does not matter which of these specific representations we choose due to the equivalence of
their first order approximations, so we will simply use the notation
$\vect{a}=\vect{a}_{(\cdot)}$ which is any one of $\vect{a}_\phi$, $\vect{a}_g$, or
$\vect{a}_p$.
\item Drop $\vect{\omega}$ from the state. This effectively ignores the effect of gyro
measurement noise relative to other error sources in the estimator - an assumption that
contributes little penalty in performance \cite{farrenkopf1978}. It also simplifies the
implementation and execution speed.
\end{enumerate}

Thus, the state is
\begin{equation}
	x = \begin{bmatrix} a \\ b\end{bmatrix}, \label{eq:state}
     	\begin{array}{l} \text{atttitude 3-vector} \\ \text{gyro bias}\end{array}
\end{equation}

The next step is to formulate expressions for the measurements. In this paper, the focus is
on vectors that are known in an earth-fixed reference frame and measured in a vehicle body
frame (by sensors mounted on the vehicle). On spacecraft, this could be star trackers or
horizon trackers, and in aircraft they are typically accelerometers and magnetometers.

The next three sections develop the EKF, CF, and UKF estimators for the state in Eq. \eqref{eq:state} assuming body-frame vector measurements.  The attitude representation results in Appendix \ref{app:att} will be used directly in this development. 

\section{Extended Kalman Filter}
Note that attitude state vector "$a$" has a singularity (as do all 3 element
representations). So we will continue to use $q$ when propagating the state vector, but the
state error covariance propagation will use $a$.

Following standard EKF formulation from \cite[Table 6.1-1]{GelbOptEst}
\begin{align}
	\dot{x}(t) &= f(x,t) + \eta(t) \label{eq:xdot}\\
	\xi_k &= h_k(x_k(t_k)) + \nu_k \label{eq:meas}
\end{align}
with
\begin{align}
	f() & \text{: The model of the time derivative of the state} \notag \\
	\eta(t) & \text{: The process noise as a function of time} \notag \\
	\xi_k & \text{: The measurement at time interval $k$} \notag \\
	h_k() & \text{: The non-linear measurement model} \notag \\
	\nu_k & \text{: The mesurement noise} \notag 
\end{align}

\subsection{State propagation}

We use Eqs. \eqref{eq:qdota}/\eqref{eq:qdotb}/\eqref{eq:qdotc} to propagate the quaternion estimate (with $\omega$ from Eq. \eqref{eq:omega}, and the convention
of denoting the estimate of a quantity with a "hat", eg. the estimate of $x$ is $\hat{x}$),
\begin{align}
	\dot{\hat{q}} = \frac{1}{2} \hat{q} \otimes \begin{bmatrix} 0 \\ \omega_g - \hat{b}\end{bmatrix}  
\end{align}
and
\begin{align}
	\dot{\hat{b}}=0 \label{eq:b_model}
\end{align}
to propagate the gyro bias estimate. To propagate the state error covariance,

\begin{equation}
f(x,t) = \begin{bmatrix} f_a \\ f_b \end{bmatrix} = \begin{bmatrix} \text{Eqs.}
\eqref{eq:aphidot},\eqref{eq:agdot},\eqref{eq:apdot} \\ 0 \end{bmatrix}
\end{equation}

\begin{equation}
	F = \frac{\partial f(x,t)}{\partial x}\bigg|_{x=\hat{x}} = 
\begin{bmatrix} \partial f_a/\partial a && \partial f_a/\partial b \\ \partial f_b/\partial
a && \partial f_b/\partial b \end{bmatrix}\bigg|_{x=\hat{x}} \label{eq:F}
\end{equation}

Now we use the first order approximation of 
\begin{equation}
	f_a = \dot{a} = \omega + \frac{1}{2} a \times \omega + \text{H.O.T.}
\end{equation}
and
\begin{equation}
	\omega \big|_{x=\hat{x}} = \omega_g - \hat{b} \label{eq:omega}
\end{equation}
where $\omega_g$ is the current gyroscope measurement vector.  Thus,
\begin{align}
\frac{\partial f_a}{\partial a} \bigg|_{x=\hat{x} }&= -\frac{1}{2} \bigl[(\omega_g -
\hat{b}) \times \bigr] \label{eq:dfaa}\\
	\frac{\partial f_a}{\partial b} \bigg|_{x=\hat{x}} &= -I \label{eq:dfab}
\end{align}
Substituting Eqs. \eqref{eq:dfaa} and \eqref{eq:dfab} into \eqref{eq:F}
\begin{equation}
F = \begin{bmatrix} -\frac{1}{2} \bigl[(\omega_g - \hat{b}) \times \bigr] && -I \\ 0 &&
0 \end{bmatrix}
\end{equation}
Now, with the state error covariance matrix (positive semi-definite, symmetric)
\begin{equation}
P = E[(x-\hat{x})(x-\hat{x})^T] = \begin{bmatrix} P_{aa} &&
P_{ab}\\P_{ba}&&P_{bb}\end{bmatrix} \label{eq:P}
\end{equation}
and $P_{aa}=P_{aa}^T, P_{bb}=P_{bb}^T,  P_{ab}=P_{ba}^T$
 we formulate the state error covariance update
\begin{align}
	\dot{P} &= FP+PF^T + Q \notag \\
	&= \begin{bmatrix}{} -(P_{ab}P_{ab}^T)  && -\frac{1}{2} [\omega \times ] P_{ab} - P_{bb}  \\
		\frac{1}{2} P_{ab}^T [\omega \times ]  - P_{bb} && 0 \end{bmatrix} + Q \label{eq:covup}
\end{align}
using $\omega$ from Eq. \eqref{eq:omega}, and
\begin{equation}
	Q = E[\eta \eta^T]
\end{equation}

\subsection{Measurement updates}

Assume a single vector measurement\footnote{It is straightforward and common to extend this to multiple simultaneous vector measurements by simply stacking the measurements into a single column vector.  We will not do this here to keep simpler notation.} in the body frame so that in Eq. \eqref{eq:meas} (dropping the measurement index $k$)
\begin{align}
	\xi &= h(q) + \nu \notag \\ 
	&= v_B(q) + \nu \label{eq:hv}
\end{align}
The subscript $B$ denotes the vector is "coordinatized" or "measured" in the $B$-frame.  Here it is helpful to be very explicit about the coordinate frames used, so we use the notation $^BC^N$ to denote the coordinate transformation of a vector in frame $N$ to the same vector in $B$.  Thus, we use the expression
\begin{align}
	v_B(q) &= ^B\negthickspace C^N(q) v_N \\
		&= ^B\negthickspace C^{\hat{B}}\negthinspace(\delta a)  \, ^{\hat{B}}C^N\negthinspace(\hat{q}) v_N \label{eq:vb} 
\end{align}

This says that the vector that is known in the $N$-frame (the earth-fixed "inertial" frame) can be transformed into the same vector in the $B$-frame by two successive transformations.  First, we transform into the current estimate of the $B$-frame which is defined by the current attitude estimate $\hat{q}$, and then we transform to the true $B$-frame using a perturbation transformation defined by the attitude perturbation (the error in the current estimate) $\delta a$.

Here the transformations are from Eqs. \eqref{eq:quattrans} and \eqref{eq:perttrans}
\begin{align}
	^{\hat{B}}C^N\negthinspace(\hat{q}) &= \bigl[ ^N C^{\hat{B}}\negthinspace(\hat{q}) \bigr]^T = C(\hat{q})^T \label{eq:CNbhat}\\
	^B\negmedspace C^{\hat{B}}\negthinspace(\delta a) &= C(\delta a)^T  = I - [\delta a \times] \label{eq:pertda}
\end{align}
Let us also define
\begin{align}
	v_{\hat{B}} &=  ^{\hat{B}}\negthickspace C^N\negthinspace(\hat{q}) v_N \label{eq:vbhat}
\end{align}

The measurement equation is then formed by substituting Eqs. \eqref{eq:vb}-\eqref{eq:vbhat} into \eqref{eq:hv} forming 
\begin{align}
	\xi &=  (I - [\delta a \times])  v_{\hat{B}} + \nu \\
	&= v_{\hat{B}} - \delta a \times v_{\hat{B}} + \nu 
\end{align}
Finally, the linearized measurement equation is
\begin{align}
	\delta \xi = \xi - v_{\hat{B}} &= v_{\hat{B}} \times \delta a + \nu  \label{eq:delxi}
\end{align}
Then the standard matrix form is
\begin{align}
	\delta \xi &= H \delta x + \nu  \notag \\
		 &= [H_a H_b] \begin{bmatrix} \delta a \\ \delta b \end{bmatrix} + \nu \label{eq:measeq}
\end{align}
with
\begin{align}
	H_a &= \frac{\partial h(x)}{\partial a} = [v_{\hat{B}} \times] \label{eq:Ha}\\ 
	H_b &= \frac{\partial h(x)}{\partial b} = 0 \label{eq:Hb}
\end{align}

Now we can form the measurement update gain matrix (again from \cite{GelbOptEst}) using Eqs. \eqref{eq:measeq}-\eqref{eq:Hb} and \eqref{eq:P}
\begin{align}
	K &= PH^T(HPH^T+R)^{-1} \\
		&= \begin{bmatrix}P_{aa}\\P_{ab}^T\end{bmatrix}H_a^T\begin{bmatrix}H_a P_{aa} H_a^T + R\end{bmatrix}^{-1} \label{eq:kalmangain}
\end{align}
where
\begin{equation}
	R = E[\nu \nu^T]
\end{equation}
so
\begin{align}
	\begin{bmatrix} \delta a \\ \delta b \end{bmatrix} &= K \delta \xi \label{eq:measupdate}\\
	\delta P &= -K H_a [P_{aa} P_{ab}]
\end{align}
then updates of covariance and bias are simply
\begin{align}
	\hat{b} &\mathrel{+}= \delta b  \label{eq:bupdate}\\
	\hat{P} &\mathrel{+}= \delta P \label{eq:Pupdate}
\end{align}
and the quaternion update from Eq. \eqref{eq:deltaq} is
\begin{align}
	\hat{q} \mathrel{+}= \frac{1}{2} \hat{q}_{last} \otimes \begin{bmatrix}0\\ \delta a\end{bmatrix} \label{eq:qupdate} 
\end{align}
If your measurement accuracy (and measurement model accuracy) warrants use of a higher order measurement update, you can use
\begin{align}
	\hat{q} = \hat{q}_{last} \otimes q(\delta a) 
\end{align}
with $q(\delta a)$ from Eq. \eqref{eq:qdela2}.

\subsection{Real measurements} \label{ssec:realmeas}

So far, we have assumed that the vector measurements in the body ($B$) frame only differ from the known vectors in the inertial ($N$) frame by a rotation plus white gaussian noise.  This may be fairly accurate for spacecraft and star trackers, but this model can have significant issues for accelerometers and magnetometers.  In this section, we discuss some of the practical approaches to help mitigate errors due to the inaccuracy of our measurement model.  

\mbox{}
\subsubsection{Accelerometer} \label{section:realmeas_accel}
First, for accelerometer measurements, a model that more accurately represents the measurements is\footnote{This model still ignores axis misalignment, scale factor, and bias errors.  These are usually present in real devices, but can often be calibrated prior to operation.  Also any offset between the point on the rigid body where the accelerometers are mounted and the point that defines $\ddot{r}$ could be accounted for, but is not shown here.}

\begin{align}
	\xi_{\text{accel}} = ^B\negthickspace C^N \negthickspace(q) (g+\ddot{r}) + \nu_{\text{accel}} \label{eq:meas_accel}
\end{align}
where the subscript $k$ denoting the measurement index is dropped.  It is clear here that the accelerometer measurement vector could differ in magnitude significantly from the gravity vector due to the translation term $\ddot{r}$.  In some cases, an estimate for $\ddot{r}$ may be available (from a position estimator that may use GPS/GNSS, for example), or a joint estimator of both translation and attitude could be employed given sufficient measurements for observability.  But it is almost always a requirement that the attitude estimator function correctly even if this enhancement is not available.  Here, we let 
\begin{align}
	v_{\hat{B}} = ^{\hat{B}}\negthickspace C^N\negthinspace(\hat{q}) (g+\ddot{r}) \label{eq:vbhat_accel}
\end{align}
Now without any normalization, the measurement update magnitude will be a function of the (signed) magnitude difference between $\xi_{\text{accel}}$ and $v_{\hat{B}}$ rather than just the angle between the two vectors, which is undesirable.  Thus we normalize both vectors to unit magnitude (represented with an overbar) and follow \eqref{eq:delxi} 
\begin{align}
	\delta \xi_{\text{accel}} = \bar{\xi}_{\text{accel}} - \bar{v}_{\hat{B}} &= \bar{v}_{\hat{B}} \times \delta a + \nu_{\text{accel}}  
\end{align}
and 
\begin{align}
	H_{a,\text{accel}} = [\bar{v}_{\hat{B}} \times] 
\end{align}
which are the accelerometer measurement forms of Eqs. \eqref{eq:delxi} and \eqref{eq:Ha}.

Admittedly, this normalization plus the errors in $\ddot{r}$ impact the gaussian assumption and the magnitude of the measurement covariance estimate $R_{\text{accel}}$.  In practice, we will set this empirically to optimize the performance of the estimator in real use-cases rather than use  noise measurements of the static accelerometer sensor only. 

\mbox{}
\subsubsection{Magnetometer} \label{section:real_mag}
A magnetometer measures the magnetic field at the sensor, and like an accelerometer, it produces a 3-element vector measured in the body frame.  The magnetic field of the earth varies depending on location and it changes slowly, but it is well modeled.  Unfortunately, the actual magnetic field can be perturbed significantly (in both magnitude and direction) from the modeled field by locally generated fields in the environment of a vehicle and by fields generated onboard the vehicle as well.  Fields generated onboard can often be calibrated and removed as long as they are static in time.  Total field variations from an earth model can be very large and can appear nearly instantaneously (if a vehicle were to drive or fly into an area influence by a local field, for example).  As a result, it is important to carefully consider how to safely include magnetometer measurements in an attitude estimator.

Here we will make use of a practical argument that will govern use of magnetometer sensing.  Due to the impact of magnetic field disturbances, we only want to use the magnetometer to aid estimation of heading (rotation about the vertical axis in the $N$ frame).  Since the gravity vector is vertical in the $N$ frame, accelerometer sensing cannot be used to sense heading.  Therefore, we need the magnetometer for heading information, but wish to isolate the magnetometer from impacting pitch and roll estimation which is available from the accelerometer.\footnote{If efforts are made to isolate magnetic disturbances and/or if you are operating in an area where such disturbances are known to be small, then of course it is possible to use the magnetic field vector to estimate any rotation perpendicular to the field direction.}

The magnetometer measurement is expressed as a function of the magnetic field vector in the $N$-frame $m$ as
\begin{align}
	\xi_{\text{mag}} = ^B\negthickspace C^N \negthickspace(q) m + \nu_{\text{mag}} \label{eq:meas_mag}
\end{align}
Now we make use of the fact that the component of a vector $m$ perpendicular to a unit vector $\mu$ can be written as
\begin{align}
	C_\mu m
\end{align}
where $C_\mu$ is the idempotent matrix
\begin{align}
	C_\mu = I - \mu \mu^T  = -[\mu\times]^2 \label{eq:cmu}
\end{align}
The second equality is from Eq. \eqref{eq:mcross2}. For our purposes here with the magnetic field measurements, we will use $\mu = [0\, 0\, 1]^T$, so $C_\mu$ is simply
\begin{align}
	C_\mu = \begin{bmatrix} 1& 0& 0 \\ 0& 1& 0 \\ 0& 0& 0\end{bmatrix}
\end{align}
but will continue the development using the general form.  The desired form can be obtained by pre-multiplying Eq. \eqref{eq:meas_mag} by 
$(^{\hat{B}}C^N) C_\mu (^NC^B)$, a matrix that transforms a vector in $B$ to $N$, zeros the component along $\mu$, then transforms back to the (estimated) body frame.
\begin{align}
	 (^{\hat{B}}C^N) C_\mu (^NC^{\hat{B}})(^{\hat{B}}C^B) \xi_{\text{mag}} = 
	 (^{\hat{B}}C^N) C_\mu  m + \tilde{\nu}_{\text{mag}} \label{eq:meas_mag2}
\end{align}
where
\begin{align}
	\tilde{\nu}_{\text{mag}}=  (^{\hat{B}}C^N) C_\mu (^NC^B)\nu_{\text{mag}} \label{eq:numag}
\end{align}
Substituting \eqref{eq:pertda} into \eqref{eq:meas_mag2} and nomalizing the vectors in the same way as with the accelerometer measurements, we obtain
\begin{align}
	\delta \xi_{\text{mag}}&= (^{\hat{B}}C^N) C_\mu [  ^N\negthickspace C^{\hat{B}} \, \bar{\xi}_{\text{mag}} - \bar{m} ] \\
	&= (^{\hat{B}} C^N) C_\mu \, ^N\negthinspace C^{\hat{B}} [\bar{\xi}_{\text{mag}} \times ] \delta a+ \tilde{\nu}_{\text{mag}} 
\end{align}
and
\begin{align}
	H_{a,\text{mag}} =^{\hat{B}}\negmedspace C^N C_\mu \, ^N\negthinspace C^{\hat{B}} [\bar{\xi}_{\text{mag}} \times ]
\end{align}
which are the magnetometer measurement forms of Eqs. \eqref{eq:delxi} and \eqref{eq:Ha}.  Again, like with the accelerometer measurements, the noise covariance is affected by the scaling and, here, by the transformation of Eq. \eqref{eq:numag}, so the covariance estimate $R_{\text{mag}} = E[\tilde{\nu}_{\text{mag}} \tilde{\nu}_{\text{mag}}^T]$ will be optimized empirically in practice.

\section{Complementary Filters}

The conceptual attractiveness of the EKF is that it promises to weight the measurements "optimally" over time.  Thus, if the state error is different in different directions, the Kalman gain will appropriately weight new measurements based on this error $P(t)$, the error estimates in the measurements $R(t)$, and the "geometry" of the measurement relative to the states being estimated $H(t)$.  However, there is simpler approach based on a geometrical observation.  Given a measurement of a vector $v_B$ in frame $B$ (denoted $\xi$ including noise) that is the known vector $v_N$ in $N$, I can transform the known vector to my current best estimate of B, $v_{\hat{B}}$, using the attitude estimate $\hat{q}$.  Then $v_B$ and $v_{\hat{B}}$ should only differ by a small rotation that can be estimated as
\begin{align} 
	\delta \theta = \bar{\xi} \times \bar{v}_{\hat{B}} \label{eq:deltatheta}
\end{align}

Now we know that this attitude error could be due to errors in the attitude estimate propagation or due to the gyro bias, but the assumption is that the gyro bias changes very slowly.  Thus we can simply update our state vector using
\begin{align}
	\delta a &= w_a \delta \theta \label{eq:deltaa} \\
	\delta b &= -w_b \delta a \label{eq:deltab}
\end{align}
where the weights $w_a$ and $w_b$ are picked empirically, and $w_b$ is a small enough to account for the slow moving assumption of this state.  Multiple vector measurements can be accommodated by combining them linearly
\begin{align}
	\delta a = \sum_i w_{a,i} \delta \theta_i
\end{align} 
The weights $w_{a,i}, i=1..n$ and $w_b$ can vary with time given any information that may be available to place confidence on respective measurements or elapsed time.  Then, the measurement update is completed as before using Eqs. \eqref{eq:bupdate} and \eqref{eq:qupdate}.  No attempt is made to estimate the state error covariance.

State propagation is the same as the EKF where Eqs. \eqref{eq:qdota}/\eqref{eq:qdotb}/\eqref{eq:qdotc} are used to propagate the quaternion estimate, 
\begin{align}
	\dot{\hat{q}} = \frac{1}{2} \hat{q} \otimes \begin{bmatrix} 0 \\ \omega_g - \hat{b}\end{bmatrix}  
\end{align}
and
\begin{align}
	\dot{\hat{b}}=0
\end{align}
to propagate the gyro bias estimate.

\subsection{Derivation starting from simplified EKF}
To help demonstrate that the EKF and this complementary filter approach have similarities, assume that the covariance matrices in the EKF have reached the following simplified form in steady state
\begin{align}
	P_{aa} &= \phi I \notag\\ P_{ab}&=-\psi I \label{eq:diagCovs}\\ R&=\rho I \notag
\end{align} 
where $\phi$, $\psi$, and $\rho$ are positive scalars.  Substituting these into Eq. \eqref{eq:kalmangain} using Eqs. \eqref{eq:Ha} and \eqref{eq:mcross2} and abbreviating $v_{\hat{B}}$ as simply $v$ for now (with $|v|^2=1$),
\begin{align}
	K &= \begin{bmatrix}\phi I\\-\psi I\end{bmatrix}[v \times]^T\biggl[\phi[v \times] I[v \times]^T + \rho I\biggr]^{-1} \\
	&=\begin{bmatrix}I\\-(\psi/\phi) I\end{bmatrix}[v \times]^T\biggl[\Bigl(1 + (\rho/\phi)\Bigr) I - v v^T\biggr]^{-1} 
\end{align} 
Then using the rank-1 version of the Matrix Inversion Lemma to compute the inverse\footnote{$(A-vv^T)^{-1}=A^{-1} + \frac{A^{-1}vv^TA^{-1}}{1-v^T A^{-1}v}, (1-v^T A^{-1}v) \ne 0$}
\begin{align}
	K &= \frac{1}{1+\rho/\phi}\begin{bmatrix}I\\-(\psi/\phi) I\end{bmatrix}[v \times]^T \biggl[I+\frac{vv^T}{\rho/\phi}\biggr] 
\end{align} 
The second term vanishes because $[v \times]^T vv^T = 0$, so
\begin{align}
	K &= \frac{1}{1+\rho/\phi}\begin{bmatrix}I\\-(\psi/\phi) I\end{bmatrix}[v \times]^T
\end{align} 
Now the measurement update is as before from Eq. \eqref{eq:measupdate} and \eqref{eq:delxi}
\begin{align}
	\begin{bmatrix} \delta a \\ \delta b \end{bmatrix} = K \delta \xi &= K (\xi - v) \\
	&= \frac{1}{1+\rho/\phi}\begin{bmatrix}I\\-(\psi/\phi) I\end{bmatrix} (\xi \times v)
\end{align}
Comparison with Eqs. \eqref{eq:deltatheta}-\eqref{eq:deltab} shows that for this example
\begin{align}
	w_a &= \frac{1}{1+\rho/\phi}\\
	w_b &= \psi/\phi
\end{align}
These results make intuitive sense in that as the measurement error becomes small relative to the attitude state error $\rho/\phi \rightarrow 0$ the weight applied to angular offset $w_a \rightarrow 1$, and conversely, as $\rho/\phi \rightarrow \infty$, $w_a \rightarrow 0$.  Further, the weight of the bias update depends on the ratio of the bias state error to the angular state error.

No claim is being made that Eqs. \eqref{eq:diagCovs} are the only way to relate the EKF to a complementary filter.  Neverthless, it shows that simple scalar approximations for covariances in the EKF do result in a complementary filter of the same form as was generated using geometric arguments.

\subsection{Stability properties of the complementary filter}
Stability of the attitude filter is of great concern when the filter output is being used for real-time dynamic control of a vehicle or robot.  This author is not aware of any stability guarantees for the EKF, however \cite{nonlinear_comp_filt_mahony_08} shows locally asymptotic stability guarantees for the complementary filter given that measurements are available from at least 2 non-colinear vectors, and \cite{jensen2011generalized} extends this to include the complementary filter with time varying positive definite matrix weights (as long as their first and second derivatives are smooth and bounded).  

When selecting an algorithm to deploy that could impact safety of operation, such stability assurances should be considered and weighed relative to performance benefits that may be available in an algorithm not guaranteed to be stable.

\section{Unscented Kalman Filter}
The UKF is related to the EKF except that it estimates the state covariance and state-measurement cross-covariance during the estimation process and does not rely on gradients that are assumed to be known.  The theoretical development and justification for the UKF can be found in \cite{julier1997new}.  

\subsection{State propagation}
We start with the current best state estimate, $\hat{q}$ and $\hat{b}$.  Then we compute a set of "sigma points" as a function of the current state error covariance that are near the current estimate, and propagate these according to the system dynamic model, then use these propagated points to estimate the new state and the new state error covariance.  Like the EKF, we perform these calculations using the 6 element state of Eq. \eqref{eq:state}, but use the dynamic model equations in the quaternion representation.  This follows the approach shown in \cite{crassidis2003unscented}.  The notation here uses an underbar to denote propagated values, and $i$ to index the sigma points (with $i=0$ for the current estimate).

First, compute the state sigma points from the columns of the positive and negative matrix square root,
\begin{align}
	\delta x(i) &= \begin{bmatrix} \delta a(i) \\ \delta b(i)\end{bmatrix} \notag \\ &= \bigl( \pm \sqrt{(n+\lambda)(P+\Delta t \cdot Q})\bigr)_i, i=1, \dotsc 2n \label{eq:sigmapoints}
\end{align}
where $\Delta t$ is the time since the last propagation and $n=6$ is the number of states in the estimator.  The parameter $\lambda$ affects the "spread" of the sigma points and can be adjusted to optimize performance.  The matrix square root ($M=\sqrt{A}$ means $M M^T=A$) can be computed with a Cholesky decomposition, and the $i$ subscript indicates the $i$th column of the matrix.  Since the state perturbation at the current state estimate is zero by definition, $\delta x(0) = 0$.

The bias propagation is simply
\begin{align}
	\delta \underline{b}(i) = \delta b(i), i=0,  \dotsc 2n \label{eq:biasprop}
\end{align}
due to the constant bias model Eq. \eqref{eq:b_model}.  The attitude quaternion propagation is three steps.  First, the quaternion sigma points are computed from Eq. \eqref{eq:sigmapoints} and the first order approximation Eq. \eqref{eq:qdela1},
\begin{align}
	q(i) = q(0) \otimes \begin{bmatrix} 1 \\ \delta a(i)/2 \end{bmatrix}, i=1,  \dotsc 2n \label{eq:quatsig}
\end{align}
with $q(0)=\hat{q}$. Next, the quaternion points are propagated (again to first order) using
\begin{align}
	\underline{q}(i) = q(i) \otimes \begin{bmatrix} 1 \\ \Delta t \cdot \omega(i)/2 \end{bmatrix}, i=0,  \dotsc 2n \label{eq:prop_q_points}
\end{align}
where the angular velocities of the sigma points are computed using the bias state sigma points (Eqs. \eqref{eq:sigmapoints} and \eqref{eq:omega}),
\begin{align}
	\omega(i) = \omega_g - (\hat{b} + \delta b(i))
\end{align}

\begin{figure}[h]
\centering
\includegraphics[scale=0.8]{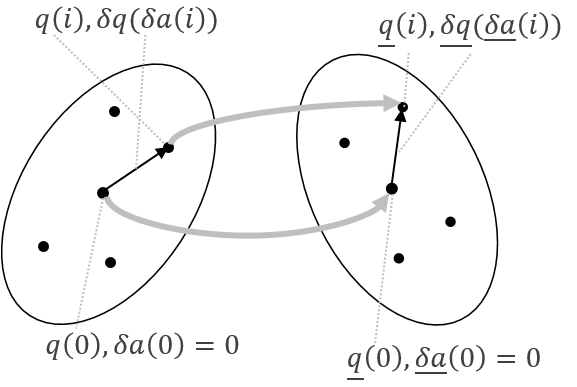}
\caption{Conceptual diagram of the UKF state propagation of quaternion states.}
\label{fig:map}
\end{figure}

The post-propagation attitude perturbation sigma points are then found from these quaternions using a form modeled after Eq. \eqref{eq:quatsig},
\begin{align}
	\underline{q}(i) = \underline{q}(0) \otimes \begin{bmatrix} 1 \\ \delta \underline{a}(i)/2 \end{bmatrix}, i=1,  \dotsc 2n \label{eq:quatbarsig}
\end{align}
Pre-multiplying this expression by the conjugate of $\underline{q}(0)$, we arrive at
\begin{align}
	\begin{bmatrix} 1 \\ \delta \underline{a}(i)/2 \end{bmatrix} = \underline{q}(0)^* \otimes \underline{q}(i), i=1,  \dotsc 2n \label{eq:quatbarsig}
\end{align}
thus $\delta \underline{a}(i)$ can be computed as twice the bottom three elements of the resulting quaternion.  (The top element should be very close to 1.0 assuming the sigma points are close to the quaternion estimate).  A conceptual representation of the quaternion sigma point propagation is shown in Fig. \ref{fig:map}.

Now the propagated sigma points are
\begin{align}
	\delta \underline{x}(i) = \begin{bmatrix} \delta \underline{a}(i) \\ \delta \underline{b}(i)\end{bmatrix} \label{eq:prop_state_sigma}
\end{align}
from Eqs. \eqref{eq:quatbarsig} and \eqref{eq:biasprop}.

Finally, the propagated state estimate and covariance estimate is computed from the propagated sigma points
\begin{align}
	\delta \hat{x} = \biggl(\frac{1}{n+\lambda}\biggr) \biggl[ \frac{1}{2} \sum^{2n}_{i=1} \delta \underline{x}(i) \biggr] \label{eq:stateperturb}
\end{align}
and the updated covariance is then
\begin{align}
	P = \biggl(\frac{1}{n+\lambda}\biggr) 
	\biggl[ \lambda \delta \hat{x} \delta \hat{x}^T + \frac{1}{2} \sum^{2n}_{i=1} (\delta \underline{x}(i) - \delta \hat{x})(\delta \underline{x}(i) - \delta \hat{x})^T \biggr]
\end{align}

To finalize the state propagation, we split the state perturbation in Eq. \eqref{eq:stateperturb} back into the attitude ($\delta \hat{a}$) and gyro bias parts ($\delta \hat{b}$) and update the full states like Eqs. \eqref{eq:qupdate} and \eqref{eq:bupdate} for the EKF.

\subsection{Measurement updates}

The vector measurement is the same as for the EKF in Eqs. \eqref{eq:hv} and \eqref{eq:vb}
\begin{align}
	\xi &= ^B\negthickspace C^N(q) v_N + \nu \label{eq:ukf_meas}
\end{align}
Now we compute the measurement sigma-points from the propagated state sigma-points of Eq. \eqref{eq:prop_q_points} (from the most recent state propagation step) 
\begin{align}
	\underline{\xi}(i) &= ^B\negthickspace C^N\bigl(\underline{q}(i)\bigr) v_N + \nu, i=0,  \dotsc 2n  \label{eq:ukf_meas2}
\end{align}
Now the measurement mean is
\begin{align}
	\hat{\xi} &= 	\biggl(\frac{1}{n+\lambda}\biggr) \biggl[ \lambda \underline{\xi}(0) + \frac{1}{2} \sum^{2n}_{i=1} \underline{\xi}(i) \biggr] \label{eq:measmean}
\end{align}
and the measurement covariance and state-measurement cross-covariance are 
\begin{align}
	P^{\xi\xi} = \biggl(\frac{1}{n+\lambda}\biggr) 
	\biggl[ \lambda (\underline{\xi}(0)-\hat{\xi})(\underline{\xi}(0)-\hat{\xi})^T + \notag
	\\ \frac{1}{2} \sum^{2n}_{i=1} (\underline{\xi}(i)-\hat{\xi})(\underline{\xi}(i)-\hat{\xi})^T  \biggr]
\end{align}
and 
\begin{align}
	P^{x\xi} = \biggl(\frac{1}{n+\lambda}\biggr) 
	\biggl[ \lambda (\delta \underline{x}(0) - \delta \hat{x})(\underline{\xi}(0)-\hat{\xi})^T + \notag
	\\ \frac{1}{2} \sum^{2n}_{i=1} (\delta \underline{x}(i) - \delta \hat{x})(\underline{\xi}(i)-\hat{\xi})^T  \biggr] \label{eq:ukfcrosscov}
\end{align}
with $\delta \underline{x}(i)$ and $\delta \hat{x}$ from Eq. \eqref{eq:prop_state_sigma} and \eqref{eq:stateperturb} respectively.

Finally, the measurement update is completed using the Kalman filter updates
\begin{align}
	K &= P^{x\xi}(P^{\xi\xi} + R)^{-1} \label{eq:ukfkalmangain}\\
	\begin{bmatrix} \delta a \\ \delta b \end{bmatrix} &= K (\xi - \hat{\xi}) \label{eq:measupdateukf}\\
	\delta P &= -K (P^{\xi\xi} + R) K^T
\end{align}
then updates of the gyro bias, covariance, and quaternion estimates are from Eqs. \eqref{eq:bupdate}, \eqref{eq:Pupdate}, and \eqref{eq:qupdate}.

\subsection{Real measurements}

Similar to Section \ref{ssec:realmeas} in the EKF development, here we add some specifics for the accelerometer and magnetometer sensors. As before, if there is a good estimate of $\ddot{r}$, the $v_N$ in Eqs. \eqref{eq:ukf_meas} and \eqref{eq:ukf_meas2} can be written as $g+\ddot{r}$.  Also, we use the same rationalle to normalize the measurement and measurement estimates.  Thus, Eq. \eqref{eq:measupdateukf} becomes
\begin{align}
	\begin{bmatrix} \delta a \\ \delta b \end{bmatrix} &= K  \delta \xi_{\text{accel}} \label{eq:measupdateukf2}\\
	\delta \xi_{\text{accel}} &= \bar{\xi}_{\text{accel}} - \bar{\hat{\xi}}_{\text{accel}}
\end{align}
with the overbar denoting the vector scaled to unit norm.  $\bar{\hat{\xi}}_{\text{accel}}$ is the normalized vector from Eq. \eqref{eq:measmean}.

For the magnetometer measurement, $v_N$  in Eqs. \eqref{eq:ukf_meas} and \eqref{eq:ukf_meas2} is replaced with $C_\mu v_N$ (using Eq. \eqref{eq:cmu}).  Then as before we can formulate the measurement 
\begin{align}
	\delta \xi_{\text{mag}} &= \bar{\xi}_{\text{mag}} - \bar{\hat{\xi}}_{\text{mag}}
\end{align}

\section{Simulations} \label{section:simulations}

To compare performance of the three different attitude estimator types, simulations were performed for a range of \textit{test cases} and \textit{scenarios} - explained in detail in the sections below.  The simulated vehicle is a 1.2 kg quadrotor with 12 Newton maximum thrust per propeller.  The consumer-grade sensors are a 3-axis accelerometer, a 3-axis gyroscope, and a 3-axis magnetometer all fixed to the body frame.  Axis alignments and scale factors are all assumed perfect (or perfectly calibrated) for these simulations.  Each accelerometer axis measurement is corrupted with 0.5 m/s$^2$ standard deviation white gaussian noise.  Each gyroscope axis is corrupted with 0.05 $^\circ$/s standard deviation white gaussian noise plus a random walk bias term.  These gyro bias terms are each generated as a gaussian random walk with 0.2 $^\circ$/s/minute standard deviation drift, and then passed through a 5 second time-constant single-pole IIR low-pass filter.  Each magnetometer axis is corrupted with 0.015 Gauss standard deviation white gaussian noise.

\noindent
\begin{table}
	\caption{Simulation test cases}
	\label{table:testcases}
	\begin{center}
		\begin{tabularx}{250pt}{| l |  X |}
		\hline
		\textbf{Test case} &  \textbf{Description} \\ \hline
		mockup\_long\_hover &  	
			No translation or rotation.  Sensors output variations are due to sensor erros only.  This case is the only case where the first and second half of the test is the same.   \\ \hline
		mockup\_easy &  
			No translation.  Angular velocities are time varying in all three axes, but are under 5 deg/sec. \\ \hline
		mockup\_slowrot &  
			No translation.  Angular velocity is fixed at 6 deg/sec in roll only.\\ \hline
		mockup &  
			No translation.  Large sinusoidal angular velocities (up to 300 deg/sec) are generated about all three axes.  \\ \hline
		straightup &  
			Vehicle is flown straight up under full throttle to altitude of just over 400 meters.  Then is it flown down as quickly as possible (without losing stability) back to ground level. \\ \hline
		bumpy\_hover &  
			Vehicle is flown with large dynamic angular velocities (up to 150 deg/sec) and translation accelerations (up to 3g), but keeping the vehicle in roughly the same position in the sky.  \\ \hline
		straightflight &  
			Vehicle is flown is a straight line at a relatively constant altitude for 700 meters, then the direction is quickly reversed and the vehicle is flown back to the starting point.\\ \hline
		longturn\_bothways &  
			Vehicle is flown into a tight coordinated turn for 1 revolutions, then the direction of the turn is reversed for 1 revolution before turning to the starting point.\\ \hline
		longturn &  
		 			Vehicle is flown into a tight coordinated turn for 2 revolutions before turning to the starting point.\\ \hline
		\end{tabularx}
	\end{center}
\end{table}

\subsection{Test cases} \label{section:testcases}
For our purposes here, we define a test case as a vehicle flight that lasts 120 seconds.  In addition to real flights we define four "mockup" cases that are not real flights, but rather a specific crafted time sequence of angular velocity and attitude without any translation accelerations.  These cases can help to validate performance for the case where the accelerometer measures the gravity vector directly (with sensor noise) without the corruption of translation accelerations.  For all test cases, the "flight" consists of a first portion (60 seconds approximately) of manuevering followed by hovering.  The second part is useful to judge the ability and speed of the of the estimator to converge to an attitude estimate of  higher accuracy than is attainable during dynamic motions of the first portion of the test.  The ground truth data (for all but the mockup cases) was produced by manually flying a quadrotor UAV ($\sim$ 1 kg total mass) in a flight training simulator and recording the resulting time sequence of vehicle position, velocity, acceleration, attitude and angular velocity.  This method of producing ground truth allows us to study highly dynamic maneuvers that are challenging for the attitude estimator, yet are realistic for this UAV vehicle.  Table \ref{table:testcases} describes the nine (9) test cases in more detail.  The descriptions are relevant for the first half of each test case with the second half being a steady hover for all cases.  The cases are ordered in the table (and later in the performance metric tables) roughly in the order increasing estimation difficulty.  For additional clarity on the character of these test cases, refer to these video clips \cite{attitude_est_gt_vids}.

\noindent
\begin{table}[t]
	\caption{Simulation scenarios}
	\label{table:scenarios}
	\begin{center}
		\begin{tabularx}{250pt}{| l |  X |}
		\hline
		\textbf{Scenario} &  \textbf{Description} \\ \hline
		nominal & 
			\begin{itemize}
				\item mag\_field\_knowledge = declination\_only  
				\item mag\_meas\_usage = horizontal   
				\item dynamic\_gains = off 
			\end{itemize}
		   \\ \hline
		no\_mag &  
			\begin{itemize}
				\item mag\_field\_knowledge = \textit{not applicable}  
				\item mag\_meas\_usage = none   
				\item dynamic\_gains = off 
			\end{itemize}
			\\ \hline
		3D\_mag &  
			\begin{itemize}
				\item mag\_field\_knowledge = xyz  
				\item mag\_meas\_usage = xyz   
				\item dynamic\_gains = off 
			\end{itemize}
			\\ \hline
		dynamic\_gains &  
			\begin{itemize}
				\item mag\_field\_knowledge = declination\_only  
				\item mag\_meas\_usage = horizontal   
				\item dynamic\_gains = on 
			\end{itemize}
			\\ \hline
		\end{tabularx}
	\end{center}
\end{table}

\subsection{Scenarios}
A \textit{scenario} is defined for our purposes here as a configuration of the estimator that is generic across all three estimator types (CF, EKF, UKF).  Specifically, the scenario is defined by a selection of each of the following:
\begin{itemize}
\item mag\_field\_knowledge = \{declination\_only, xyz\}.  This determines the estimator's knowledge of the local earth magnetic field.  This is typically found using an onboard lookup table from the vehicle's location on earth, to either a declination\footnote{The angular separation between true and magnetic north in the horizontal plane.} (``declination\_only") or a full three dimensional field estimate (``xyz'').  Such information can be obtained using a model such as the World Magnetic Model \cite{world_mag_model}.  Of course these models have limited resolution resulting in some local errors, and would not include any magnetic anomalies or locally generated fields.  If a network connection is available on the vehicle, the information can be obtained over the network or the locally stored tables can be updated in real-time.
\item mag\_meas\_usage = \{none, horizontal, xyz\}.  This determines the use of the magnetic measurement in the estimator.  ``None'' means that magnetic measurements are ignored.  This would result in unbounded drift of the heading estimate without additional sensing sources.  ``horizontal'' means that only the field in the local horizontal plane will be used, such as the technique described in Section \ref{section:real_mag} 
\item dynamic\_gains = \{off, on\}.  Since the accelerometer measures the gravity vector much more accurately when the vehicle translation accelerations are small, a technique that weights the accelerometer measurements with the nominal weights during high dynamics but with higher weights during low dynamic portions of a flight is attractive.  ``off'' means the nominal case of a single static set of weights.  ``on'', for the purposes of this paper, means that a second set of weights can be employed during low dynamic motions.  The detector of low vs. high dynamics is very important for this to work well (else this type of technique can end up being detrimental to overall performance).  The detector used here is outlined in Section \ref{section:dyndetect}.  The concept of weights are different for the different estimator types.  For CF, these are simply a second set of weights, but for EKF and UKF the change is in the measurement error covariance matrix $R$ of Eqs. \eqref{eq:kalmangain} and \eqref{eq:ukfkalmangain} respectively. 
\end{itemize}

Among the possible permutations of the items above, we have chosen to define the four scenarios described in Table \ref{table:scenarios} for simulation.

\subsection{Simulation parameters}
\noindent
\begin{table}
	\caption{Simulation parameters}
	\label{table:simparams}
	\centering
	\begin{center}
		\begin{tabular}{| m{0.4cm} | m{0.8cm} | m{0.6cm} | c | m{0.6cm} |} \hline
		\textbf{Est} &  \multicolumn{3}{c|}{\textbf{Parameters}} & \textbf{Eqs.}\\ \hline
		\multirow{4}{0.5cm}{CF}  &  \multirow{4}{0.8cm}{Meas weights} & accel & 
					$w_a=0.0002\begin{bmatrix} 1\\1\\1\end{bmatrix}$$ w_b=0.03\begin{bmatrix} 1\\1\\1\end{bmatrix}$ & \multirow{3}{.7cm}{\eqref{eq:deltaa}, \eqref{eq:deltab} }  \\ \cline{3-4}
		                     &   & low-dyn accel & $w_a=0.002\begin{bmatrix} 1\\1\\1\end{bmatrix}$$ w_b=0.03\begin{bmatrix} 1\\1\\1\end{bmatrix}$ &   \\ \cline{3-4}
		                     &   & mag & $w_a=0.002\begin{bmatrix} 1\\1\\1\end{bmatrix}$$ w_b=0.03\begin{bmatrix} 1\\1\\1\end{bmatrix}$ &   \\ \cline{3-4}
		                     &   & low-dyn mag  & $w_a=0.02\begin{bmatrix} 1\\1\\1\end{bmatrix}$$ w_b=0.03\begin{bmatrix} 1\\1\\1\end{bmatrix}$ &   \\ \hline
		\multirow{5}{*}{EKF}  &  \multirow{5}{0.8cm}{Meas cov} & accel & 
					$R=\text{diag}\begin{bmatrix} \alpha^2 \\\alpha^2 \\\alpha^2 
                                                                               \end{bmatrix}, \begin{matrix} \alpha=30 \text{ m/s$^2$}\end{matrix}$
				& \multirow{4}{0.7cm}{\eqref{eq:kalmangain} }  \\ \cline{3-4}
		                     &   & low-dyn accel & 
					$R=\text{diag}\begin{bmatrix} \alpha^2 \\\alpha^2 \\\alpha^2 
                                                                               \end{bmatrix}, \begin{matrix} \alpha=4.9 \text{ m/s$^2$}\end{matrix}$
					 &   \\ \cline{3-4}
		                     &   & mag & 
					$R=\text{diag}\begin{bmatrix} \alpha^2 \\\alpha^2 \\\alpha^2 
                                                                               \end{bmatrix}, \begin{matrix} \alpha=0.47 \text{ gauss}\end{matrix}$
					 &   \\ \cline{3-4}
		                     &   & low-dyn mag & 
					$R=\text{diag}\begin{bmatrix} \alpha^2 \\\alpha^2 \\\alpha^2 
                                                                               \end{bmatrix}, \begin{matrix} \alpha=0.094 \text{ gauss}\end{matrix}$
					 &   \\ \cline{2-5}
		                     &  Process noise cov & \multicolumn{2}{c|}{
					$Q=\text{diag}\begin{bmatrix} \alpha^2 \\\alpha^2 \\\alpha^2 \\
						                                         \beta^2\\\beta^2\\\beta^2 
                                                                               \end{bmatrix}, \begin{matrix} \alpha=1 \text{ $^\circ$/s} \\ \beta=0.5 \text{ $^\circ$/s/min}\end{matrix}$
					} & \eqref{eq:covup} \\ \hline
		\multirow{6}{*}{UKF}  &  \multirow{6}{1cm}{Meas cov} & accel &
					$R=\text{diag}\begin{bmatrix} \alpha^2 \\\alpha^2 \\\alpha^2 
                                                                               \end{bmatrix}, \begin{matrix} \alpha=10 \text{ m/s$^2$}\end{matrix}$
				& \multirow{5}{0.7cm}{\eqref{eq:ukfkalmangain} }  \\ \cline{3-4}
		                     &   & low-dyn accel & 
					$R=\text{diag}\begin{bmatrix} \alpha^2 \\\alpha^2 \\\alpha^2 
                                                                               \end{bmatrix}, \begin{matrix} \alpha=4.9 \text{ m/s$^2$}\end{matrix}$
					 &   \\ \cline{3-4}
		                     &   & mag & 
					$R=\text{diag}\begin{bmatrix} \alpha^2 \\\alpha^2 \\\alpha^2 
                                                                               \end{bmatrix}, \begin{matrix} \alpha=0.24 \text{ gauss}\end{matrix}$
					 &   \\ \cline{3-4}
		                     &   & low-dyn mag & 
					$R=\text{diag}\begin{bmatrix} \alpha^2 \\\alpha^2 \\\alpha^2 
                                                                               \end{bmatrix}, \begin{matrix} \alpha=0.094 \text{ gauss}\end{matrix}$
					 &   \\ \cline{2-5}
		                     &  Process noise cov & \multicolumn{2}{c|}{
					$Q=\text{diag}\begin{bmatrix} \alpha^2 \\\alpha^2 \\\alpha^2 \\
						                                         \beta^2\\\beta^2\\\beta^2 
                                                                               \end{bmatrix}, \begin{matrix} \alpha=1 \text{ $^\circ$/s} \\ \beta=0.5 \text{ $^\circ$/s/min}\end{matrix}$
					}  &  \eqref{eq:sigmapoints} \\ \cline{2-5}
		                     &  Overall & \multicolumn{2}{c|}{$\lambda=1$}  & \eqref{eq:sigmapoints}-\eqref{eq:ukfcrosscov}  \\ \hline
		\end{tabular}
	\end{center}
\end{table}

Estimator weights (weight vectors for CF, and process and measurement noise covariances for both EKF and UKF) were set first from the models, and then refined manually through extensive simulations to balance performance between the high-dynamics and low-dynamics portions of the test cases (see Section \ref{section:testcases} for further explanation of these cases).  For EKF and UKF, using covariances based only on the linear models is not effective during high-dynamic maneuvers where the attitude errors can grow beyond typical small-angle, linear assumption.  So manual tuning of these estimator covariances was required.  Estimator weights for CF are tuned manually in any case, but this tuning is simplified by the fact that there is only a single weight for each sensor (for both angle and bias weights).

Final weighting parameters used for the simulations are compiled in Table \ref{table:simparams}.  The "diag" function used in the table (to save space) constructs a diagonal matrix using the elements of the specified vector as the diagonal elements.  Parameters marked "low-dyn" are only used in scenarios with dynamic\_gains="on" based on the output of the dynamics detector of Section \ref{section:dyndetect}.

\subsection{Performance metrics}

The performance metrics that we are using for evaluating the results were selected for the particular test cases and scenarios used and are summarized in Table \ref{table:metrics}.  The first two are most useful for evaluating performance in the initial, high-dynamic phase of the test cases.  MaxEVz is the absolute value of the maximum value of the z-component of the Euler vector error and MaxEVxy is the absolute value of the maximum of both the x and y-components of the Euler vector error.  Here, since the vehicle can deviate significantly from level flight (pitch/roll near zero), it is problematic to try to use Euler angle errors to judge estimator accuracy (these errors are a function of the vehicle attitude).  Thus we use Euler vectors instead.  Given an attitude estimate expressed as a quaternion, $\hat{q}$ and the corresponding ground-truth quaternion, $q$, compute the error as
\begin{align}
	\delta q =
		\begin{bmatrix} \delta q_0 \\ \vect{\delta q} \end{bmatrix} &= 
		\hat{q}^* \otimes q \label{eq:quat_err}
\end{align}
Now the error Euler vector, $\delta \vect{a}_\phi$,  can be expressed by combining Eq. \eqref{eq:quat_err} with Eqs. \eqref{eq:ev} and \eqref{eq:quat} resulting in
\begin{align}
	\delta \vect{a}_\phi = \frac{\phi}{\delta q_0 \tan(\phi/2)} \vect{\delta q} \label{eq:}
\end{align}
with \begin{align}
	\phi = 2 \cos^{-1}(\delta q_0)
\end{align}
thus
\begin{align}
	\text{MaxEVz} &= \bigl| \max ( (\delta \vect{a}_\phi)_{z} ) \bigr|\\
	\text{MaxEVxy} &= \bigl| \max \bigl( \max ( (\delta \vect{a}_\phi)_{x},(\delta \vect{a}_\phi)_{y}) \bigr) \bigr|
\end{align}
Note that this metric is undefined for exactly zero attitude error, but this is not a practical issue in tests with any noise.  We split the z-axis from xy since the z performance is mostly impacted by the magnetometer when the vehicle is close to level, so this split is only useful for test cases like mockup\_long\_hover, mockup\_easy, and straightup where the dynamic portion of the test is close to level.

The second pair of metrics is useful for evaluating steady-state performance in level hover, and it also validates convergence after the dynamic portion of the flight.  FinH is the final (at the end of the test) heading error, and FinPR is the maximum of the final pitch and roll error.
\begin{align}
	\text{FinH} &= \bigl| \text{heading} (t_{\text{final}}) \bigr|\\
	\text{FinPR} &= \bigl| \max \bigl( \text{pitch}(t_{\text{final}}),\text{roll}(t_{\text{final}}) \bigr) \bigr|
\end{align}
\noindent
\begin{table}
	\caption{Performance metrics}
	\label{table:metrics}
	\begin{center}
		\begin{tabularx}{250pt}{| l |  X |}
		\hline
		\textbf{Metric} &  \textbf{Description} \\ \hline
		MaxEVz (deg) &  	
		Compute the sequence of attitude errors expressed as Euler vectors in the body frame.  Take the maximum of the z-axis components.
		   \\ \hline
		MaxEVxy  (deg) &  
		Compute the sequence of attitude errors expressed as Euler vectors in the body frame.  Take the maximum of the x and y-axis components.
			\\ \hline
		FinH (deg) &  
		The final error in heading
			\\ \hline
		FinPR (deg) & 
		The maximum of the final errors in pitch and roll
			\\ \hline
		\end{tabularx}
	\end{center}
\end{table}

\subsection{Dynamics detector}
\label{section:dyndetect}
The dynamics detector is a causal algorithm that determines when the estimator should more heavily weight the sensor information, namely when low-dynamics=TRUE.  It is a conservative detector in that it requires confidence to build over time before it will produce this output.

Let us define the current absolute deviation of the accel norm from the gravitational constant (9.81 m/s$^2$) as
\begin{align}
	\mu = \bigl| |a| - |g| \bigr|
\end{align}  
and also define $\mu_f$ as a low-pass filtered version of this value.  Given a high threshold, $T_H$, and a low threshold, $T_L$, 
\begin{align*}
	\text{if} (\mu > T_H), &\text{  reset  } \mu_f \text{  to  } \mu \\
	\text{if} (\mu_f < T_L), & \\
	                                   & \text{  low-dynamics = TRUE} \\
	\text{else}, & \\
	                                   & \text{  low-dynamics = FALSE} \\
\end{align*}  
For the simulation results shown here for scenario ``dynamic\_gains,'' this detector was used with a 5 second moving average for the low-pass filter, and with $T_H=2.0 m/s^2$ and $T_L=0.7 m/s^2$.

\subsection{Discussion of results}

The results shown here are a subset of the number of cases studied --- and is the subset that the author sees as the most interesting and instructive.  For example, all the results here assume the same sensor degradation models (as described at the beginning of this section, Section \ref{section:simulations}).  Tests with perfect inertial sensing and, alternately, with very poor inertial sensing, plus tests with magnetometer biases of up to 10 degrees were performed.  Of course, the estimator settings must change for optimal performance with different sensor models, and the resulting performance depends on the magnitude of the sensor errors, but other than this, there were no surprises.  Further, magnetometer offsets resulted in the expected offsets in attitude estimates, but again there were no surprises or particularly interesting side effects.  Thus, these results were omitted to save space.  Additionally, techniques for incorporating translation acceleration (as described in Section \ref{section:realmeas_accel}) were incorporated and simulated, but when noise was added to these estimates, using them in the estimator produced significant degradation of overall performance.  This is not to say that this technique could not be effectively integrated, but noise in these translation acceleration inputs would need to be sufficiently small, and time delays between these estimates and the inertial sensors measurements would need to be carefully managed.  We did not quantify the required noise levels or acceptable latency values for the approach to be effective, however.

\input{./estimator_results_nominal/performance_table.tex}
\input{./estimator_results_no_mag/performance_table.tex}
\input{./estimator_results_3D_mag/performance_table.tex}
\input{./estimator_results_dynamic_gains/performance_table.tex}

Note also that each result presented in the tables and plots is from a single statistical realization of all the random sources in the simulation.  Monte-carlo simulations were not seen as necessary for drawing the general conclusions targeted in this paper.  Further, there are quite a few test cases presented, so the reader should be able to form a good idea of the range of results even without being given a distribution of the final metrics from repeated simulation trials.   This is mentioned to be clear that individual metrics for a test case are not necessarily representative of a statistical quantity such as the  mean.

\begin{figure}[!ht]
\centering
\includegraphics[trim= 0mm 25mm 0mm 5mm,scale=0.6,clip]
{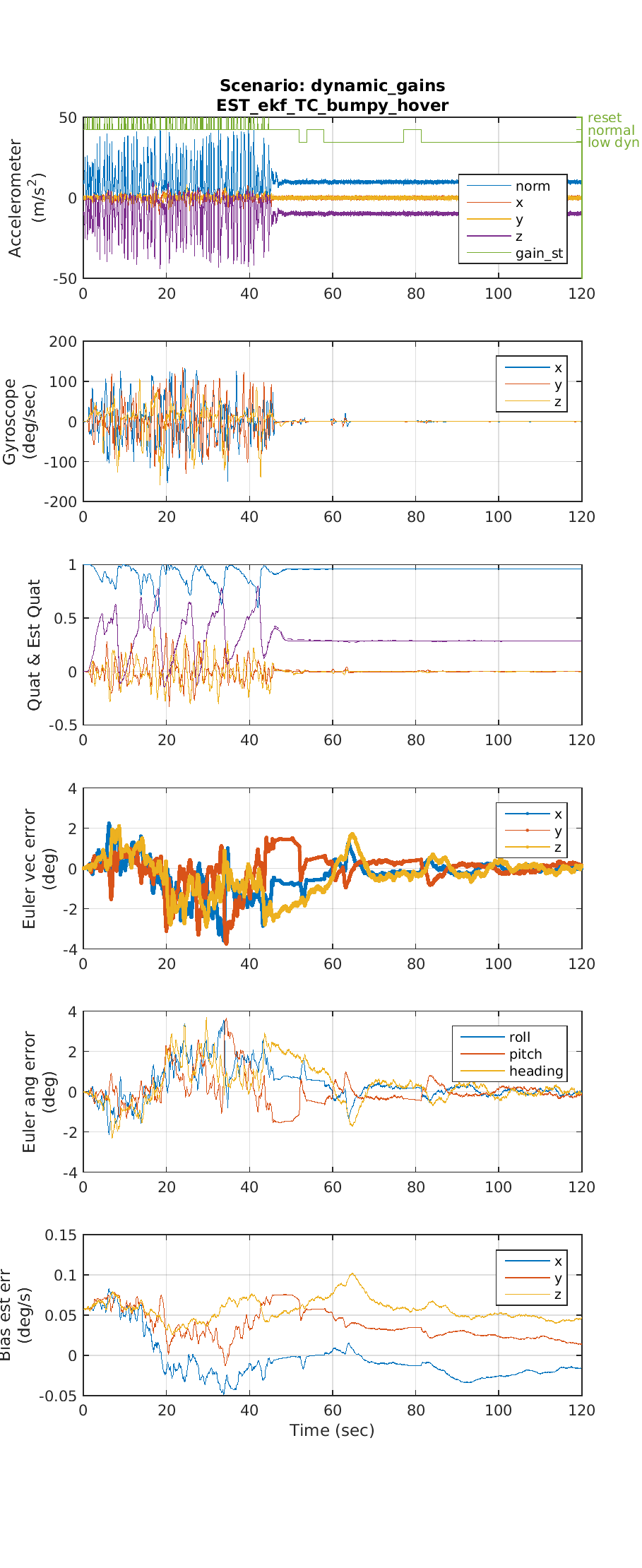}
\caption{}
\label{fig:Scen_dynamic_gains_Est_ekf_TC_bumpy_hover_SC_nominal_inertial_paper}
\vspace{20mm}
\end{figure}

The discussion below will refer to time-plots with a common format.  Referring to Fig. \ref{fig:Scen_dynamic_gains_Est_ekf_TC_bumpy_hover_SC_nominal_inertial_paper} as an example, the title specifies the scenario, the estimator type, and the test case. The first and second subplots show the raw inertial measurements from the accelerometers and gyroscopes in the body-fixed x, y, and z-axes.  For plots of the scenario "dynamic\_gains," the accel plot also shows the 3-state output of the dynamics detector of Section \ref{section:dyndetect} with vertical-axis labels on the right of the plot (and denoted in the legend as "gain\_st" for "gain state."  The third subplot shows the truth quaternion as solid lines, and the estimated quaternion (the estimated attitude as a quaternion) as dashed lines.  The fourth and fifth subplots show the angular errors of the attitude estimate as Euler vector errors, and Euler angle errors, respectively.  The sixth and final subplot shows the gyro bias estimate error for each body fixed axis. 

Now we discuss some observations from the results of the simulations.  A small number of selected plots is included here in this paper for discussion, but the full set of results is available on the internet in Ref. \cite{attitude_est_plots}.

The use of the accelerometer measurements as a reference for pitch and roll is fundamentally limited when the sensor experiences translation accelerations.  Thus during highly dynamic motion where these translation accelerations are large, the attitude estimation performance suffers for all estimator types.  In the tests shown here, a single set of estimator parameters (those shown in Table \ref{table:simparams}) were selected to balance the performance of both the dynamic portions, and the "near hover" portions of each test.  The exception to this is the dynamic\_gains scenario where two sets were used and were selected based on the output of the dynamics detector. 

Referring to the tabulated metrics in Tables \ref{table:nominal}-\ref{table:dynamic_gains}, for the EKF and CF, max absolute errors of up to 20 degrees were common for the high dynamics portions of the test cases with significant translation acceleration (namely, bumpy\_hover, straightflight, longturn\_bothways, and long\_turn).  The UKF suffered substantially worse performance during these highly dynamic periods.  After the vehicle returned to near-hover, all estimators converged to sub-degree accuracy, with a few exceptions with slightly higher residual error.

The poor performance of the UKF estimator for this setup is not fully understood, but given that the error covariance is estimated from the data, it is possible that this error covariance is poorly estimated under these conditions, thus impacting the overall UKF attitude estimator performance. 

\begin{figure}[!ht]
\centering
\includegraphics[trim= 0mm 25mm 0mm 5mm,scale=0.6,clip]
{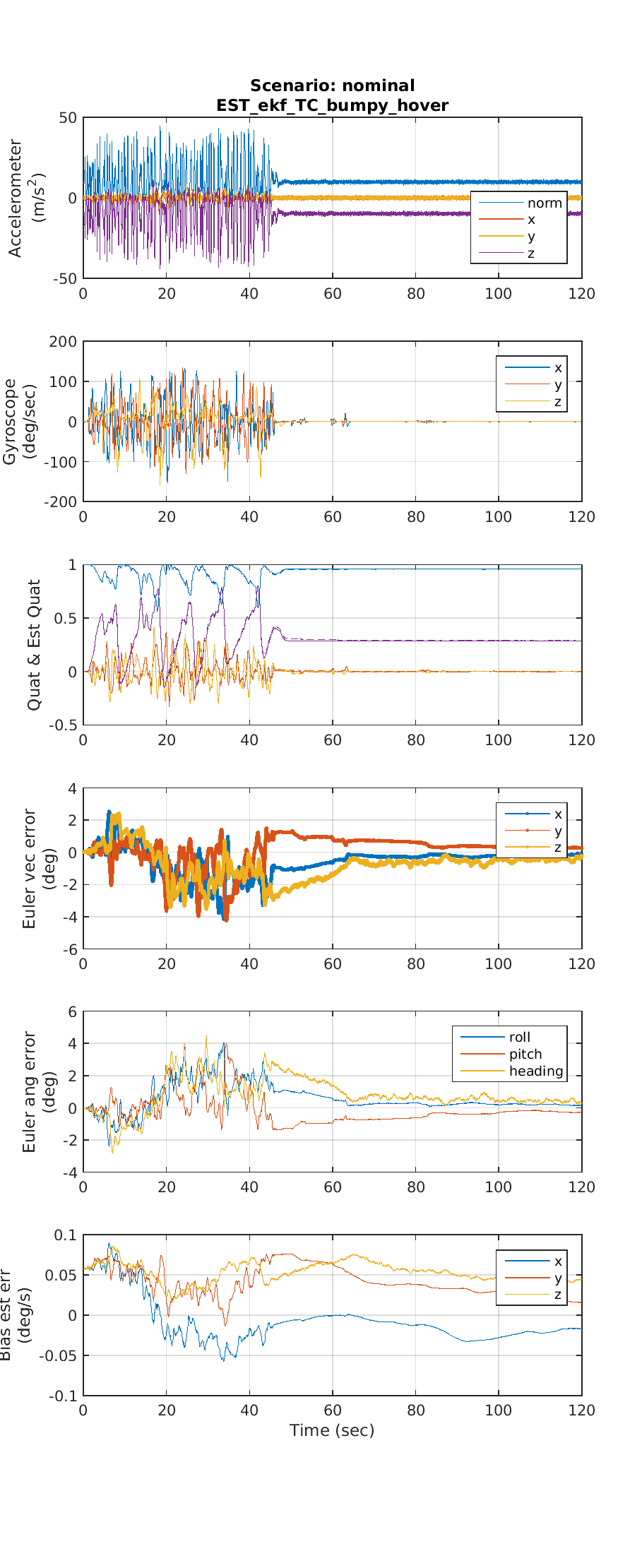}
\caption{}
\label{fig:Scen_nominal_Est_ekf_TC_bumpy_hover_SC_nominal_inertial_paper}
\vspace{20mm}
\end{figure}

\begin{figure}[!ht]
\centering
\includegraphics[trim= 0mm 25mm 0mm 5mm,scale=0.6,clip]
{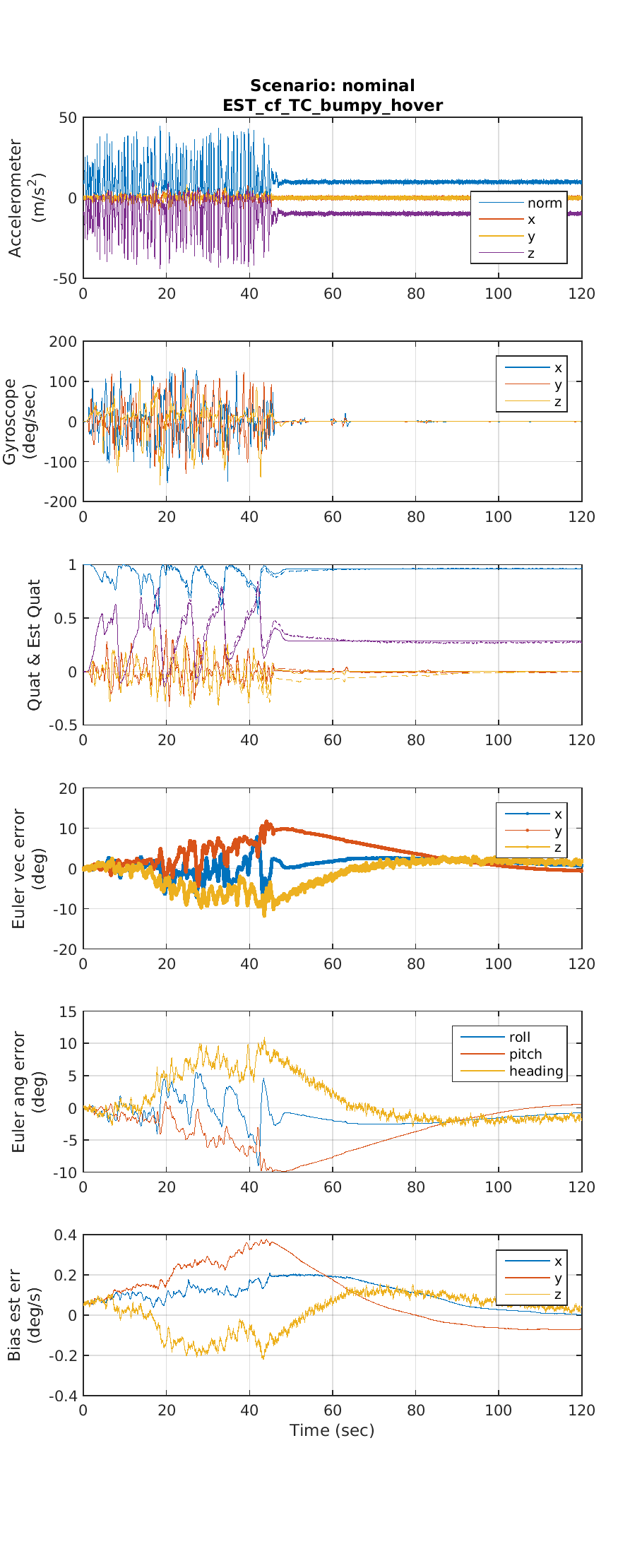}
\caption{}
\label{fig:Scen_nominal_Est_cf_TC_bumpy_hover_SC_nominal_inertial_paper}
\vspace{20mm}
\end{figure}

\begin{figure}[!ht]
\centering
\includegraphics[trim= 0mm 25mm 0mm 5mm,scale=0.6,clip]
{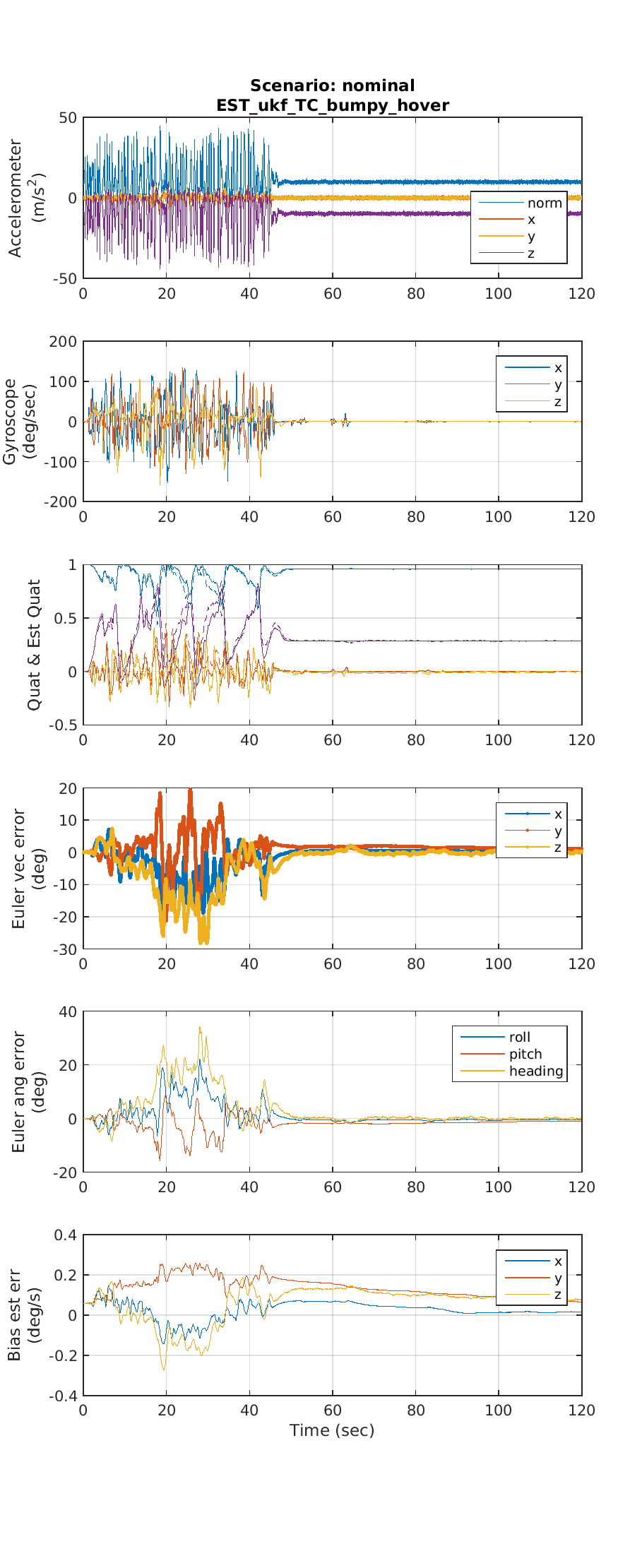}
\caption{}
\label{fig:Scen_nominal_Est_ukf_TC_bumpy_hover_SC_nominal_inertial_paper}
\vspace{20mm}
\end{figure}

\begin{figure}[!ht]
\centering
\includegraphics[trim= 0mm 25mm 0mm 5mm,scale=0.6,clip]
{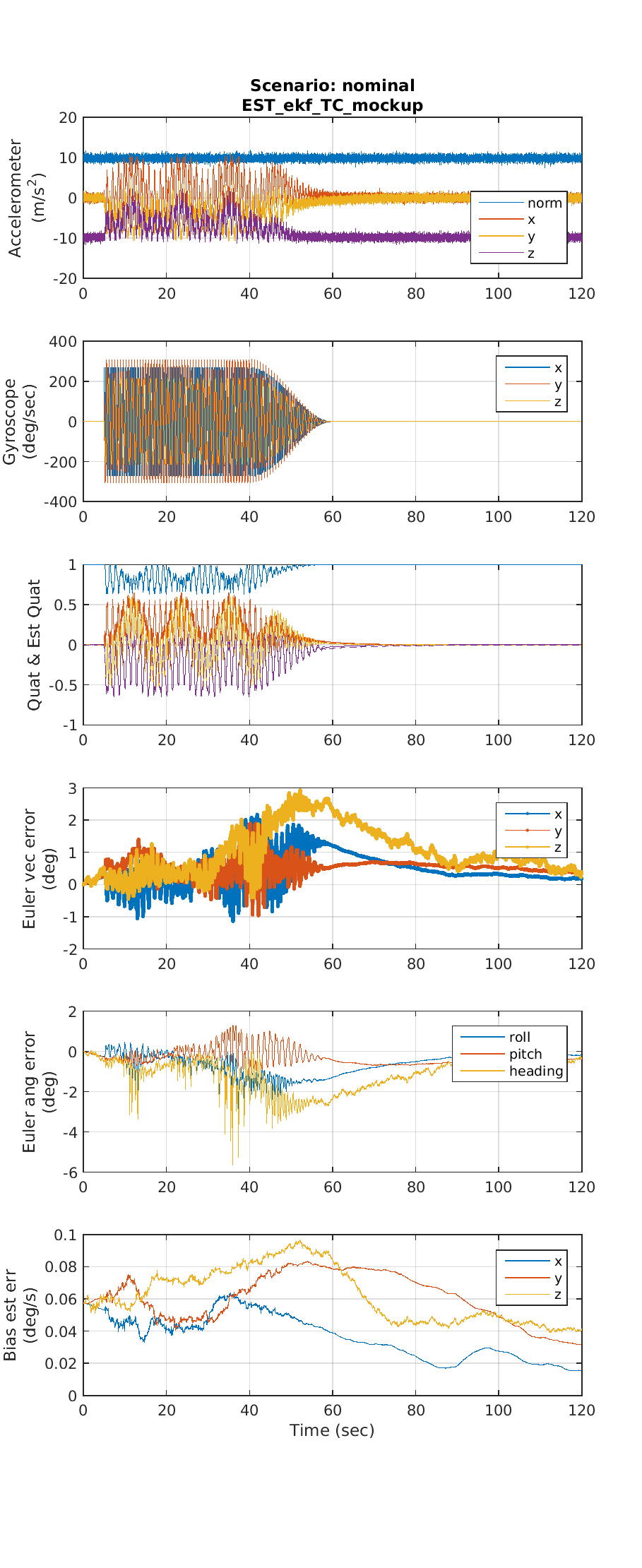}
\caption{}
\label{fig:Scen_nominal_Est_ekf_TC_mockup_SC_nominal_inertial_paper}
\vspace{20mm}
\end{figure}

\begin{figure}[!ht]
\centering
\includegraphics[trim= 0mm 25mm 0mm 5mm,scale=0.6,clip]
{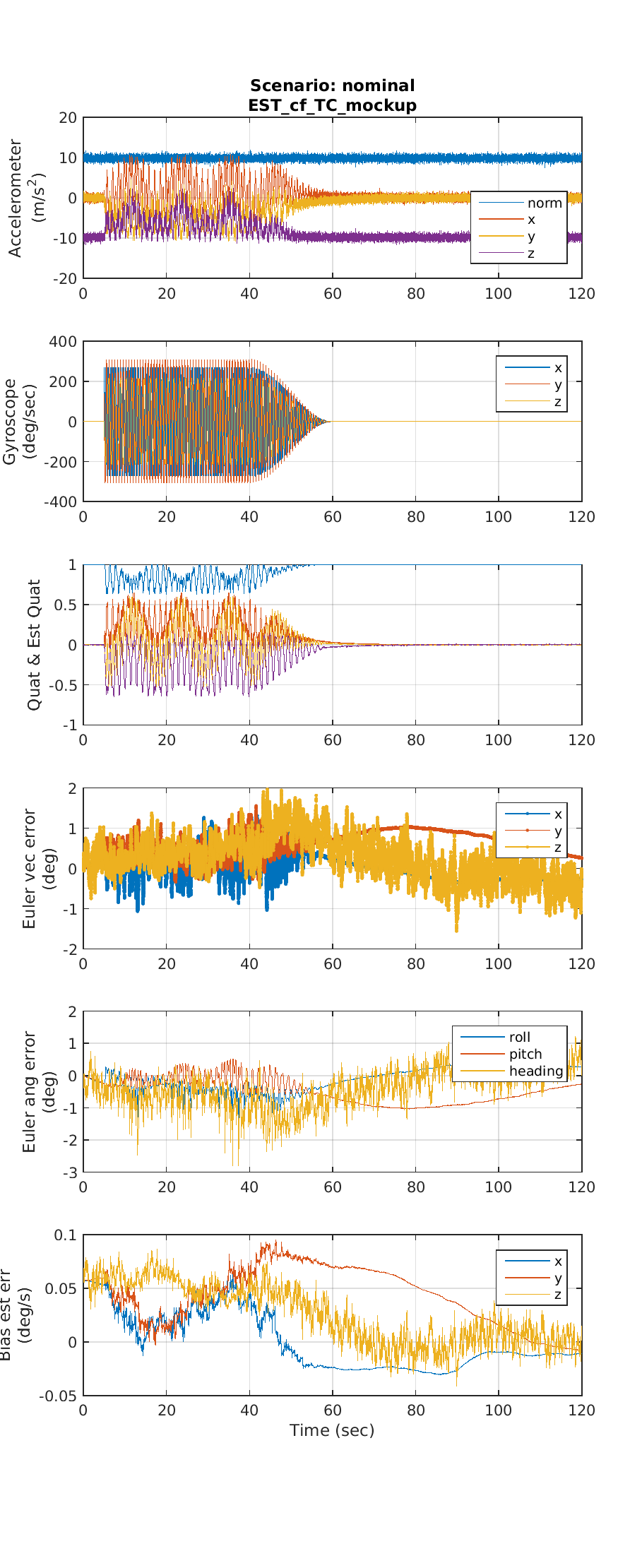}
\caption{}
\label{fig:Scen_nominal_Est_cf_TC_mockup_SC_nominal_inertial_paper}
\vspace{20mm}
\end{figure}

\begin{figure}[!ht]
\centering
\includegraphics[trim= 0mm 25mm 0mm 5mm,scale=0.6,clip]
{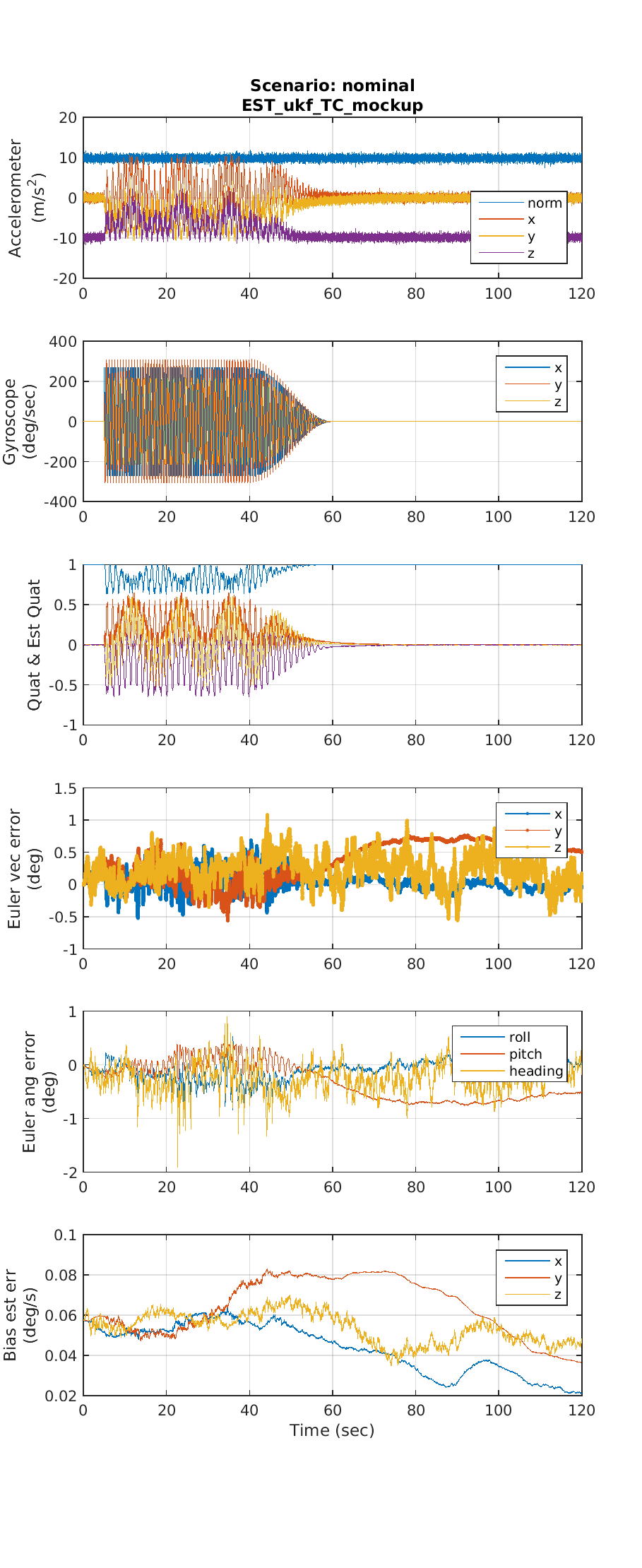}
\caption{}
\label{fig:Scen_nominal_Est_ukf_TC_mockup_SC_nominal_inertial_paper}
\vspace{20mm}
\end{figure}

\begin{figure}[!ht]
\centering
\includegraphics[trim= 0mm 25mm 0mm 5mm,scale=0.6,clip]
{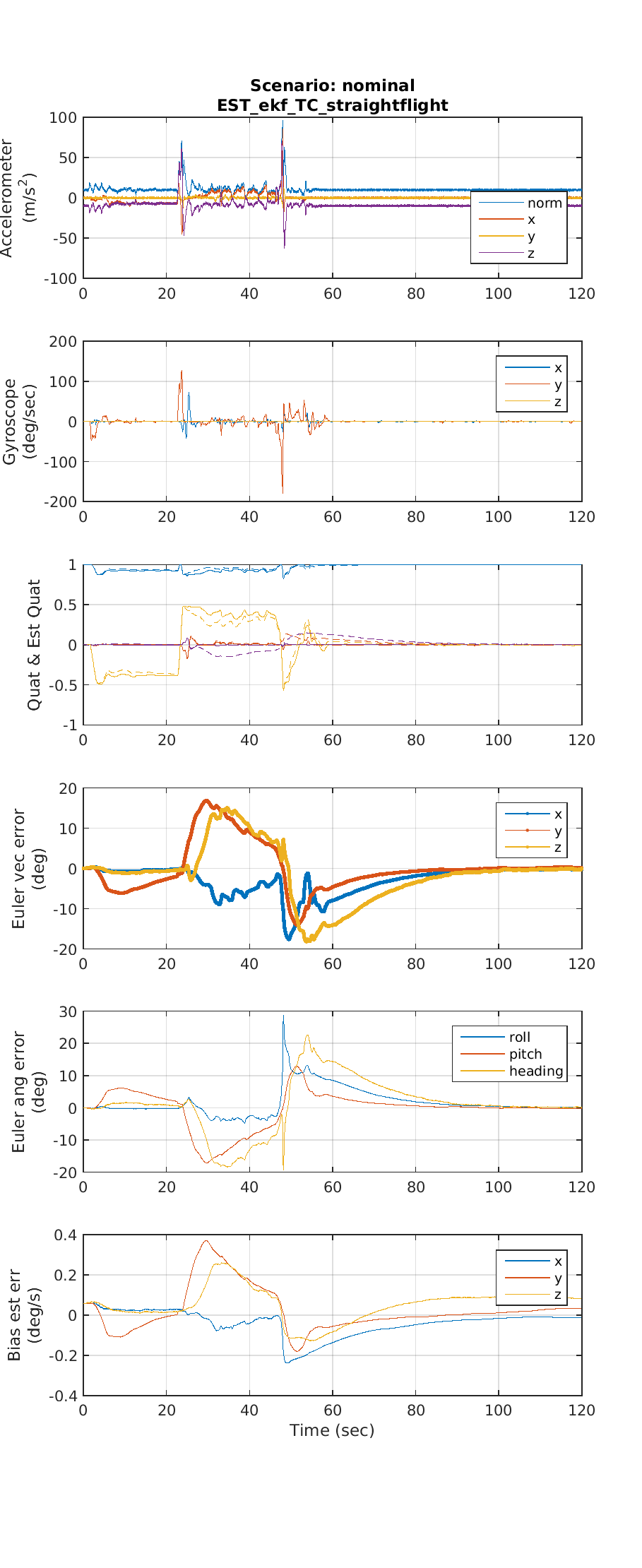}
\caption{}
\label{fig:Scen_nominal_Est_ekf_TC_straightflight_SC_nominal_inertial_paper}
\vspace{20mm}
\end{figure}

\begin{figure}[!ht]
\centering
\includegraphics[trim= 0mm 25mm 0mm 5mm,scale=0.6,clip]
{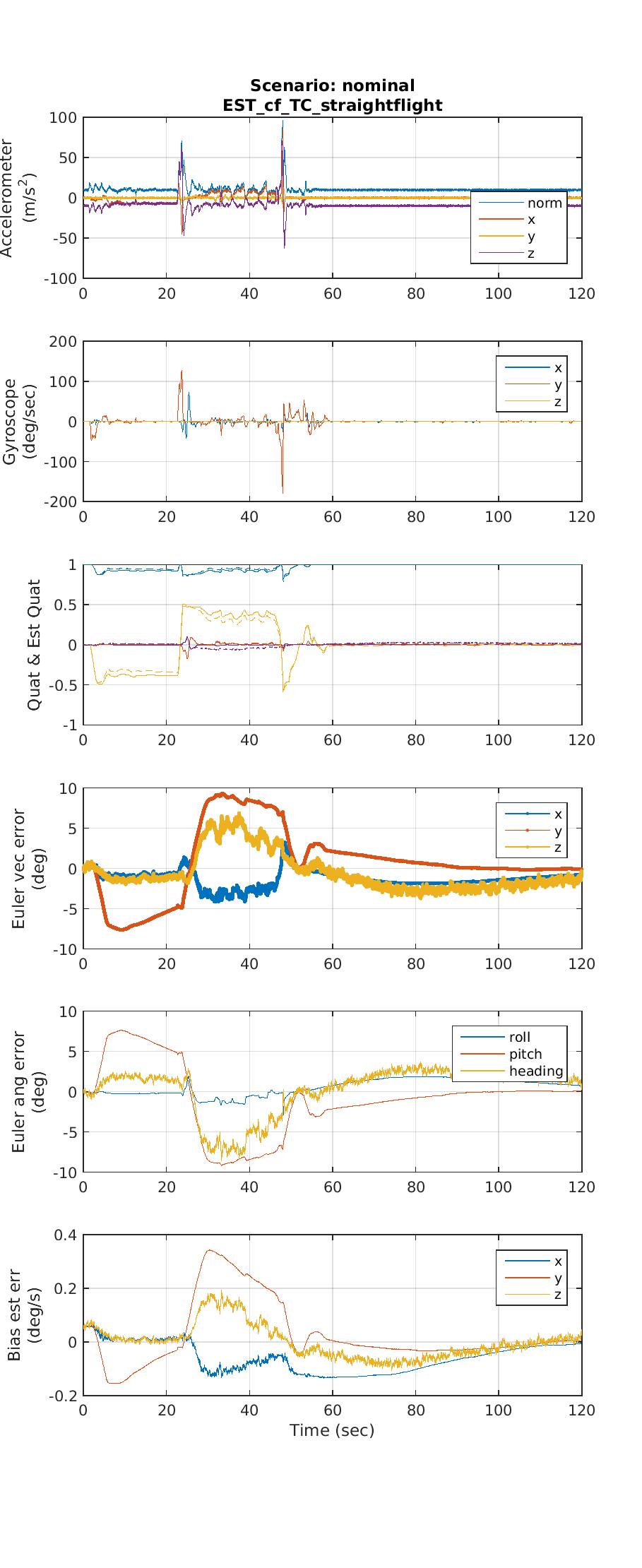}
\caption{}
\label{fig:Scen_nominal_Est_cf_TC_straightflight_SC_nominal_inertial_paper}
\vspace{20mm}
\end{figure}

\begin{figure}[!ht]
\centering
\includegraphics[trim= 0mm 25mm 0mm 5mm,scale=0.6,clip]
{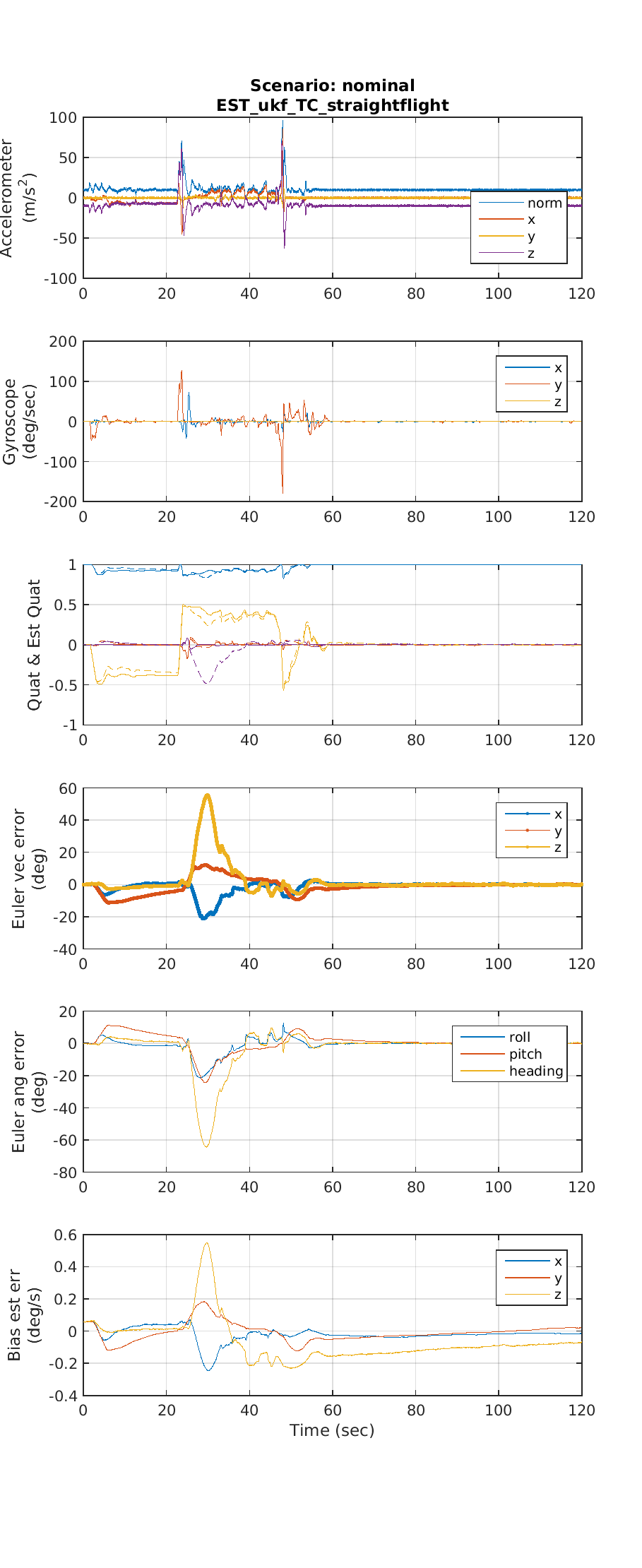}
\caption{}
\label{fig:Scen_nominal_Est_ukf_TC_straightflight_SC_nominal_inertial_paper}
\vspace{20mm}
\end{figure}

Taking a closer look at the plots, first we focus on the nominal scenario for the test cases bumpy\_hover, mockup and straightflight.  Both bumpy\_hover and mockup  involve highly dynamic attitude motions without much overall translation.  Of course the bumpy\_hover case has significant translation accelerations and is based on a real simulated flight while the mockup case has zero translational accelerations (and the angular velocities are prescribed by sinusoudal functions).  Straightflight is an out and back flight in a straight line but with significant accelerations at the start, stop, and the turnarnound point.  Figures \ref{fig:Scen_nominal_Est_ekf_TC_bumpy_hover_SC_nominal_inertial_paper} - \ref{fig:Scen_nominal_Est_ukf_TC_straightflight_SC_nominal_inertial_paper} show the full results for each of the estimator types and these are consistent with the tabulated result in Table \ref{table:nominal}.

For the bumpy hover case, the EKF is the clear winner with Euler vector errors within $\sim$4 degrees throughout while the CF errors become a bit over 10 degrees.  However, for the straightflight case, the CF constrains the maximum errors (9 degrees) better than the EKF (18 degrees) during the large accelerations.  Again, the UKF is significantly worse than the others, but all three converge similarly during hover.  For the mockup test case, all three estimators perform similarly (but with the UKF demonstrating the best overall performance), highlighting that the translational accelerations are critical for testing the limits of performance.  Comparing the EKF and the CF, these results indicate that they have comparable performance despite the increased complexity of the EKF.  While the EKF is clearly better for the bumpy\_hover test case, the CF is better for straightflight, and they are about the same for longturn\_bothways and longturn.  Note that this conclusion may not extend to other setups with better sensor performance, so other work demonstrating benefits of an EKF are not being challenged here in the general sense.

\begin{figure}[!ht]
\centering
\includegraphics[trim= 0mm 25mm 0mm 5mm,scale=0.6,clip]
{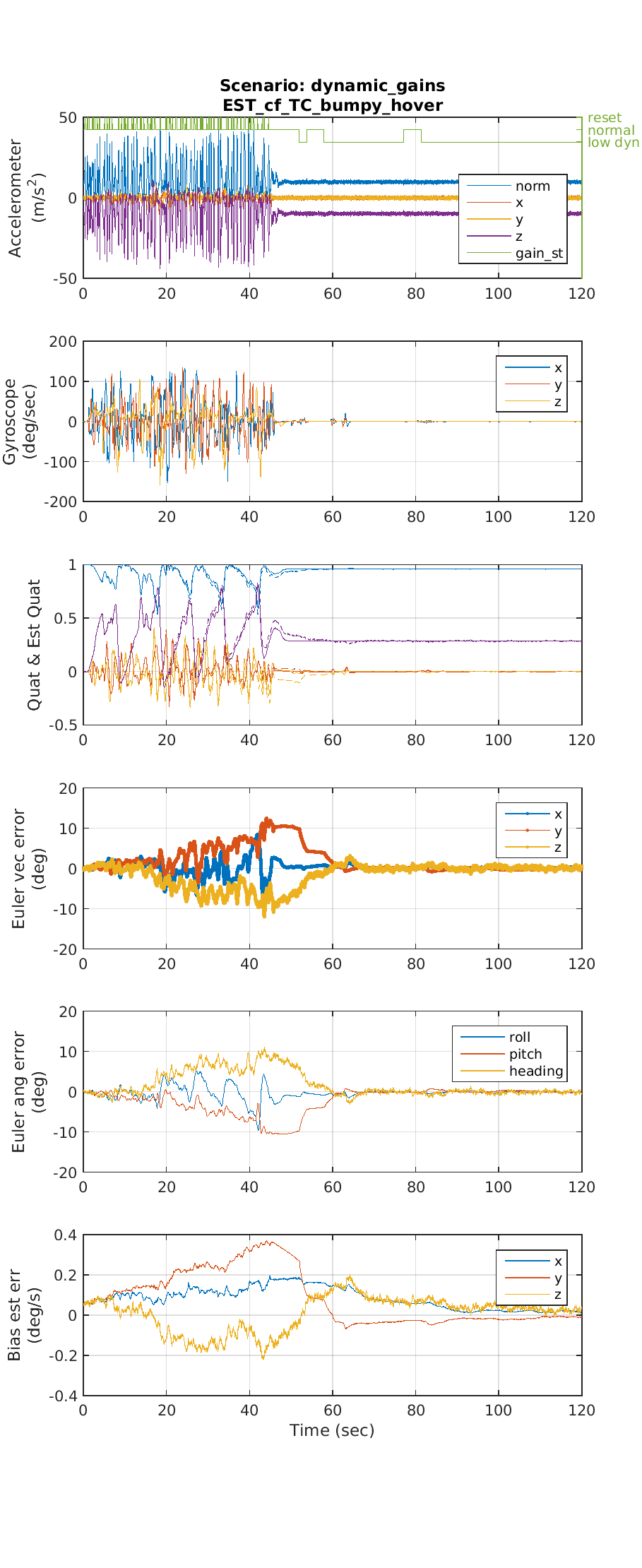}
\caption{}
\label{fig:Scen_dynamic_gains_Est_cf_TC_bumpy_hover_SC_nominal_inertial_paper}
\vspace{20mm}
\end{figure}

Next, looking at the dynamic\_gains scenario, it is clear from comparing the CF estimator for the bumpy\_hover case (Figs. \ref{fig:Scen_nominal_Est_cf_TC_bumpy_hover_SC_nominal_inertial_paper} and \ref{fig:Scen_dynamic_gains_Est_cf_TC_bumpy_hover_SC_nominal_inertial_paper}), that the error converges much more quickly when the "low dynamics" is detected and the larger measurement gains are used.  In looking at the Euler vector error from $t=50-60$ seconds on the plot, the convergence rate corresponds directly to the output of the dynamics detector output shown in the first subplot (the accelerometer plot, right-hand axis).  It is tempting to use those gains throughout due to the improved performance in hover, however this was tried and resulted in dramatically worse performance during the dynamic motions of the test cases.  Note that during development, some tuning was required to the dynamics detector to make sure that the larger gains were not applied during the high dynamics portions.  Before this tuning, the dynamic\_accel scenario produced extremely poor results, so care must be taken by an implementer to avoid this.

For the magnetometer scenarios, the tabulated results are sufficient so not plots are included here.  Comparing Tables \ref{table:nominal} and \ref{table:no_mag} shows that when the magnetometer is not used (scenario no\_mag), the heading results suffer as expected.  In fact the heading drifts in a random walk so comparing the absolute numbers is not particularly valuable.  The results do confirm that the pitch and roll results are not significantly affected by removing the magnetometer measurement (metrics MaxEVxy and FinPR).  For the 3D\_mag scenario, the magnetometer provides not only heading information but also pitch and roll information (the ambiguity is the rotation about the field vector) and the full 3-axis earth fixed field vector is assumed known at the vehicle location.  Comparing Tables \ref{table:nominal} and \ref{table:3D_mag} shows that this information can improve results during high-dynamics.  For example, looking at the straightflight test case and the MaxEVxy metrics, the benefit of 3D\_mag is clear for the EKF (17.61 deg reduces to 1.01 degrees) and CF (9.28 degrees reduces to 2.89 degrees).  Benefits are also clear for the bumpy\_hover test case.   However, the assumption for this scenario was that the magnetic field vector was known perfectly at the location of the vehicle.  Any errors in the earth mag model or from locally generated (especially man-made ) magnetic fields would compromise this benefit. 

It is worth emphasizing that this paper focuses more on comparing estimator performance for a baseline set of algorithms in challenging conditions rather than finding the overall limits of performance.  Thus, additional sensing sources such as computer vision may significantly change the overall performance of all algorithms.  Further, it is expected that a more sophisticated algorithm for adjusting filter gains/weights than the simple dynamics detector of Section \ref{section:dyndetect} could also provide significant benefits in performance.

\begin{figure}[!ht]
\centering
\includegraphics[trim= 0mm 25mm 0mm 5mm,scale=0.6,clip]
{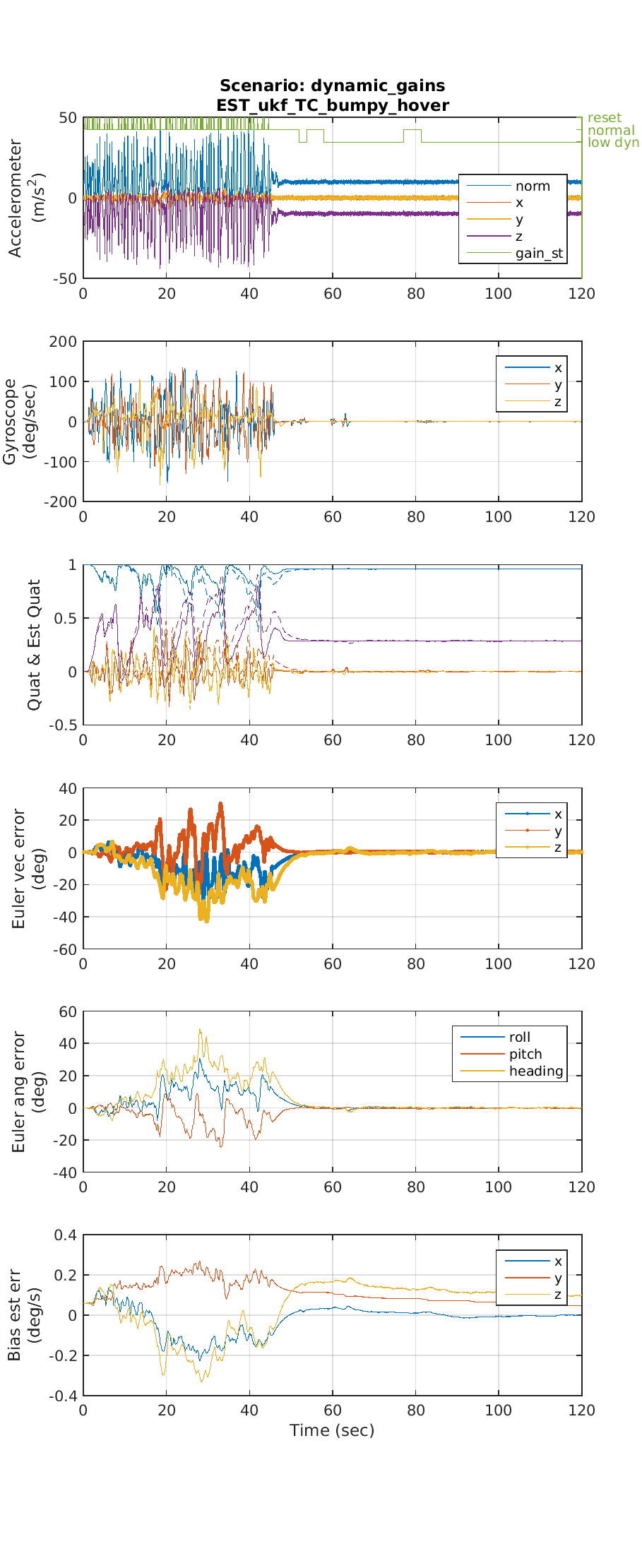}
\caption{}
\label{fig:Scen_dynamic_gains_Est_ukf_TC_bumpy_hover_SC_nominal_inertial_paper}
\vspace{20mm}
\end{figure}

\section{Conclusions}

For ``near-hover'' conditions, excellent attitude estimation performance is achieved by all three estimator types (EKF, CF, and UKF).  This is especially true when the estimators use parameters that are optimized for the near-hover case.  For example, the FinPR and FinH metrics in Table \ref{table:dynamic_gains} show $\sim\negthickspace\frac{1}{2}$ degree accuracy for all test cases and $\sim\negthickspace\frac{1}{10}$ degree for many.  However, during highly dynamic maneuvers, the translation accelerations sensed (in addition to the gravitational acceleration) degrade the attitude estimation performance significantly.  For example, maximum attitude errors between 10-20 degrees were commonly observed in the simulations for the EKF and CF.  The UKF significantly underperformed these estimators often producing 2 times the attitude errors during high dynamics.

Different magnetometer scenarios were tested in simulation, and intuitive predictions confirmed.  Namely, the elimination of a magnetometer had little to no negative impact on pitch and roll, but allowed the heading estimate to drift (random walk) to large errors.  Addditionally, using all three axes of the magnetometer measurement was very effective for the EKF and CF with the caveat that a perfect magnetic model was assumed.  Realistic magnetic model errors would need to be added to confirm an overall benefit, but this was not attempted here.

A simple dynamics detector was implemented to allow the estimators to place stronger weights on the accelerometer measurement when translation accelerations were expected to be small.  While careful tuning of this estimator was required, this approach resulted in the overall best performance for all estimator types.

\newpage
\appendices
\section{Attitude Representations}
\label{app:att}

\subsection{Definitions}
The most common attitude representations in use are direction cosine matrices, euler angles,
the euler vector, the quaternion, the Gibbs vector, and MRP or modified Rodrigues
parameters. Here we avoid Euler angles since they involve somewhat complex trigonometric caluculations in their kinematics (see the appendix of \cite{KLLSpacecraftDyn}).\footnote{They are quite useful, however, in intuitively understanding a particular attitude.  For example, there is nothing easier than pitch, roll, and heading to understand the orientation of a vehicle.  Thus, we often convert attitude to these parameters for final visulaization and plotting even though the filters and controllers do not use them.}  The direction cosine matrix is used, but since it is not a vector form, it is not as obvious how to use it in a standard Kalman filter formulation, and with 9 parameters specifying 3 degrees of freedom, it is quite redundant (through the orthonormal constraint).

Here we focus on the three 3-element representations (Euler vector, Gibbs vector, and MRP) plus the quaternion.  The quaternion is valuable since it globally non-singular, though with its 4 elements, it is constrained (to unit magnitude).  The others have a singularity (as do all 3 value attitude represenatations), but are unconstrained.  They are valuable in filters mostly for small perturbations where their singularity does not arise.  See \cite[Introduction]{markley2003attitude} for more discussion of the theoretical justifications for this choice in optimal filtering.  While it would be possible to simply focus on one of the 3-element representations, here we retain all three.  This is mainly because it is interesting to compare the geometric and kinematic expressions between the representations.  We also demonstrate that their first and second-order approximations are often equivalent, with the proper scaling.   

Using Euler's theorem that any coordinate transformation can be achieved by a single
rotation about an axis fixed in both the inital and final reference frames, we define:
\begin{align}
	\text{Euler vector}		&&	\vect{a}_\phi &= \phi \vect{e} \label{eq:ev} \\
	\text{Gibbs vector}		&&	\vect{a}_g &= 2 \tan(\phi/2) \vect{e} \label{eq:gv}\\
	\text{MRP}			&&	\vect{a}_p &= 4 \tan(\phi/4) \vect{e} \label{eq:mrp}\\
	\text{Quaternion}		&&	q =
		\begin{bmatrix} q_0 \\ \vect{q} \end{bmatrix} &= 
		\begin{bmatrix} \cos(\phi/2) \\ \sin(\phi/2) \vect{e} \end{bmatrix}  \label{eq:quat}
\end{align}

\noindent 
This rotation should be thought of as aligning a reference frame with the initial frame, and then rotating it about the axis defined by unit vector $\vect{e}$ by $\phi$ radians (using the right-hand rule for the sign of the rotation) to the final frame. For vector transformation, the rotation is in the opposite direction, thus a vector coordinatized in the final frame can be transformed to its coordinatization in the initial frame by rotating it $\phi$ radians about $\vect{e}$.

 The first three of these of these representations are [3x1] and the quaternion is [4x1] and unit norm by
definition. For the [3x1] vectors, the scale factor is chosen such that the magnitude of the
vector approaches the rotation angle for small angles. Such small angle attitudes are
represented using the attitude perturbation notation $\vect{\delta a}$, thus $ |
\vect{\delta a}_{(\cdot)} | =\phi$.

It is also useful to express the quaternion in terms of the [3x1] vectors $\vect{a}_\phi$, $\vect{a}_g$, and $\vect{a}_p$.  From Eqs. \eqref{eq:ev}-\eqref{eq:quat}, the first is trivial and the last two take some algebraic and trigonometric manipulation that is not shown here.
\begin{align}
	q(\vect{a}_\phi) &= \begin{bmatrix} \cos(\alpha_\phi) \\[3pt] \dfrac{ \vect{a}_\phi }{2\alpha_\phi} \sin(\alpha_\phi)\end{bmatrix}  \label{eq:qaphi} \\[5pt]
	q(\vect{a}_g) &= \frac{1}{\sqrt{1+\alpha_g}}\begin{bmatrix} 1 \\[3pt] \dfrac{\vect{a_g}}{2} \end{bmatrix}  \label{eq:qag} \\[5pt]
	q(\vect{a}_p) &= \frac{1}{1+\alpha_p}\begin{bmatrix} 1-\alpha_p \\[3pt] \dfrac{\vect{a_p}}{2} \end{bmatrix}  \label{eq:qap} 
\end{align}
with the following parameters defined to simplify the expressions
\begin{align}
	\alpha_\phi &= \frac{1}{2}| \vect{a_\phi}| \\[5pt]
	\alpha_g     &= \frac{1}{4} | \vect{a_g}|^2 \\[5pt]
	\alpha_p     &= \frac{1}{16} | \vect{a_p}|^2
\end{align}

Taylor series expansion of these expressions shows that all three representations have identical first and second order approximations for small rotations.  To first order
\begin{equation}
	q(\vect{\delta a}_{(\cdot)}) = \begin{bmatrix} 1 \\ \vect{\delta a}_{(\cdot)}/2 \end{bmatrix}  \label{eq:qdela1} \\
\end{equation}
and to second order
\begin{equation}
	q(\vect{\delta a}_{(\cdot)}) = \begin{bmatrix} 1-|\vect{\delta a}_{(\cdot)}|^2/8 \\ \vect{\delta a}_{(\cdot)}/2 \end{bmatrix}  \label{eq:qdela2} \\
\end{equation}

\subsection{Kinematics}

The following kinematic expressions are the exact derivatives of attitude with respect to
the current attitude and the angular velocity (of the final frame in the inital frame). The first equation is in vector form, and the
second equation is the equivalent matrix form. For more details and derivations, see
\cite{shuster1993survey}. Note also that the first order approximation for small rotations 
is the same for each of the [3x1] attitude representations
\begin{equation}
\vect{\dot{\delta a}}_{(\cdot)} = \omega + \frac{1}{2} \vect{\delta a}_{(\cdot)} \times
\vect{\omega} = \Bigl[ I + \frac{1}{2} \bigl[ \vect{\delta a}_{(\cdot)} \times \bigr] \Bigr]
\vect{\omega}
\end{equation}

\uline{Euler vector}
\begin{align}
\vect{\dot{a}}_\phi &= \omega + \frac{1}{2} \vect{a}_\phi \times \vect{\omega} +
\frac{1}{4\alpha_\phi^2} \bigl( 1-\alpha_\phi\cot\bigl(\alpha_\phi \bigr) \bigr) \vect{a}_\phi \times
(\vect{a}_\phi \times \vect{\omega)} \label{eq:aphidot}\\
&= \Bigl[ I + \frac{1}{2} \bigl[ \vect{a}_\phi \times \bigr] + \frac{1}{4\alpha_\phi^2} \bigl( 1-\alpha_\phi\cot\bigl(\alpha_\phi \bigr) \bigr)
			\bigl[ \vect{a}_\phi \vect{a}^T_\phi - I a^2 \bigr] \Bigr] \vect{\omega} 
\end{align}

\uline{Gibbs vector}
\begin{align}
\vect{\dot{a}}_g &= \omega + \frac{1}{2} \vect{a}_g \times \vect{\omega} + \frac{1}{4}
(\vect{\omega} \cdot \vect{a}_g) \vect{a}_g \label{eq:agdot}\\
&= \Bigl[ I + \frac{1}{2} \bigl[ \vect{a}_g \times \bigr] + \frac{1}{4} \bigl[ \vect{a}_g
\vect{a}^T_g \bigr] \Bigr] \vect{\omega}
\end{align}

\uline{MRP}
\begin{align}
\vect{\dot{a}}_p &= (1- \alpha_p) \omega + \frac{1}{2}
\vect{a}_p \times \vect{\omega} + \frac{1}{4} (\vect{\omega} \cdot \vect{a}_p) \vect{a}_p
\label{eq:apdot} \\
&= \biggl[ (1- \alpha_p) I + \frac{1}{2} \bigl[ \vect{a}_p
\times \bigr] + \frac{1}{4} \bigl[ \vect{a}_p \vect{a}^T_p \bigr] \biggr] \vect{\omega}
\end{align}

\uline{Quaternion} (using Eq. \eqref{eq:quatmult})
\begin{align}
	\dot{q} &= \frac{1}{2} q \otimes \begin{bmatrix} 0 \\ \vect{\omega} \end{bmatrix} \label{eq:qdota}\\
	 &=\frac{1}{2}\begin{bmatrix} -\vect{q} \cdot \vect{\omega} \\ \vect{q}\times \vect{\omega} +q_0 \vect{\omega}  \end{bmatrix} \\
	&= \frac{1}{2}\begin{bmatrix} -\vect{q}^T \\ \bigl[ \vect{q} \times \bigr] + q_0 I \end{bmatrix}\vect{\omega} \label{eq:qdotb}
\end{align}

An equivalent expression that is simple to evaluate is
\begin{align}
	\dot{q} &= \frac{1}{2} \Omega q \label{eq:qdotc}
\end{align}

with,
\begin{align}
	\Omega &= \begin{bmatrix}
			0 && -\omega_1 && -\omega_2 && -\omega_3\\
			\omega_1 && 0 && \omega_3 && -\omega_2\\
			\omega_2 && -\omega_3 && 0 && \omega_1\\
			\omega_3 && \omega_2 && -\omega_1 && 0	
			\end{bmatrix}
\end{align}

\uline{Direction Cosine Matrix}
\begin{align}
		\dot{C} = C \bigl[ \vect{\omega} \times \bigr]
\end{align}

\subsection{Direction Cosine or Transformation Matrix}

The following are the expressions for the direction cosine or transformation matrix as a
function of the attitude representation. Let $C$ be the matrix that transforms a vector
expressed in the final frame to the initial frame.\footnote{Note that this definition of the "unlabeled" direction cosine matrix is the opposite of the convention in some references, for example \cite{shuster1993survey}.  Care is required because authors who adopt that convention also may redefine the quaternion multiplication operator using a negative sign for the cross product term in Eq. \eqref{eq:quatmult} to retain the same order for cascade transormations like Eqns. \eqref{eq:cascade_q} and \eqref{eq:cascade_C}.} Also, similar to the kinematic expressions, we can show that the first order expressions of the transformation matrix is
the same for each of the [3x1] attitude representations
\begin{equation}
	C(\vect{\delta a}_{(\cdot)}) =  I  + [\vect{\delta a}_{(\cdot)}\times] \label{eq:perttrans}
\end{equation}
In fact the second order expressions for the direction cosine matrices are equivalent as well
\begin{equation}
	C(\vect{\delta a}_{(\cdot)}) =  I  + [\vect{\delta a}_{(\cdot)}\times]+\frac{1}{2}[\vect{\delta a}_{(\cdot)}\times]^2 \label{eq:perttrans2}
\end{equation}

\uline{Euler vector}
\begin{equation}
	C(\vect{a}_\phi) = I + \sin\phi[\vect{e}\times]+ (1-\cos\phi)[\vect{e}\times]^2
\end{equation}

\uline{Gibbs vector} 
\begin{equation}
C(\vect{a}_g) = I +\frac{1}{1+ \alpha_g }  [\vect{a}_g\times]+\frac{1}{2(1+ \alpha_g)}  [\vect{a}_g\times]^2
\end{equation}

\uline{MRP} 
\begin{equation}
C(\vect{a}_p) =  I + \frac{1-\alpha_p}{(1+\alpha_p)^2} [\vect{a}_p\times]+
\frac{1}{2(1+\alpha_p)^2} [\vect{a}_p\times]^2
\end{equation}

\uline{Quaternion}
\begin{align}
	C(q) &=  \bigl(q_0^2 - |\vect{q}|^2\bigr)  I  + 2q_0[\vect{q}\times]+ 2\vect{q}\vect{q}^T \label{eq:quattrans} \\
	 &=   I  + 2q_0[\vect{q}\times]+ 2[\vect{q}\times]^2
\end{align}

Note that some references in the literature differ in these expressions in scale factors and using the identity
\begin{equation}
	[\vect{v}\times]^2 = \vect{v}\vect{v}^T-|\vect{v}|^2 I \label{eq:mcross2}
\end{equation}

\subsection{Quaternion operations}

The quaternion multiplication operator is defined as
\begin{equation}
	q \otimes p = \begin{bmatrix} q_0 p_0 - \vect{q} \cdot \vect{p} \\ q_0 \vect{p}+p_0 \vect{q} + \vect{q} \times \vect{p} \end{bmatrix} \label{eq:quatmult}
\end{equation}
Cascading transformations through intermediate frame $A$ uses quaternion multiplication
\begin{align}
	^Nq^B = ^N\negthickspace q^A \otimes ^A\negthickspace q^B \label{eq:cascade_q}
\end{align}
just like cascading direction cosine matrices
\begin{align}
	^NC^B = ^N\negthickspace C^A \, ^AC^B \label{eq:cascade_C}
\end{align}
The quaternion conjugate is defined as
\begin{equation}
	q^* = \begin{bmatrix} q_0 \\ - \vect{q} \end{bmatrix}
\end{equation}
such that
\begin{equation}
	q^* \otimes q = q \otimes q^* = q_I = \begin{bmatrix} 1\\0\\0\\0 \end{bmatrix}
\end{equation}
where $q_I$ is the "identity" quaternion.  Quaterion multiplication is distributive over addition and associative, but not commutative.  Changing the order of the arguments results in
\begin{equation}
	p \otimes q = (q^* \otimes p^*)^*
\end{equation}

When the quaternion is a tranformation from, say, frame $B$ to frame $N$, denoted $^Nq^B$, then a vector coordinatized in the $B$ frame can be converted into its coordinatization in $N$ (transformed) using
\begin{equation}
	\begin{bmatrix} 0 \\ v_N \end{bmatrix} = ^N \negmedspace q^B \otimes \begin{bmatrix} 0 \\ v_B \end{bmatrix} \otimes (^Nq^B)^*  
\end{equation}
or
\begin{equation}
	\begin{bmatrix} 0 \\ v_N \end{bmatrix} = ^N \negmedspace q^B \otimes \begin{bmatrix} 0 \\ v_B \end{bmatrix} \otimes ^B\negmedspace q^N  
\end{equation}
since
\begin{equation}
	^B q^N = (^Nq^B)^*    
\end{equation}

This is equivalent to
\begin{equation}
	 v_N  = C(^Nq^B) v_B   
\end{equation}
with $C(q)$ from Eq. \eqref{eq:quattrans}.

When cascading transformations from $B$ to intermediate frame $A$ then to $N$
\begin{equation}
	\begin{bmatrix} 0 \\ v_N \end{bmatrix} = ^N\negmedspace q^A \otimes ^A\negmedspace q^B \otimes \begin{bmatrix} 0 \\ v_B \end{bmatrix} \otimes ^Bq^A \otimes ^Aq^N  
\end{equation}

It is also useful to look at changes in quaternions, since filters typically compute quaternion state updates.  Let an initial quaternion $q$ be changed by an update quaternion $q_\Delta$ forming final $\bar{q}$ 
\begin{equation}
	\bar{q} =  q \otimes  q_\Delta
\end{equation}
This can also be written as an update to the elements of $q$ we will call $\Delta q$.  It's important to note that $\Delta q$ is just 4 numbers and not a unit quaternion. 
\begin{equation}
	\bar{q} = q+\Delta q  
\end{equation}
Combining these two equations and using the distributive property
\begin{align}
	\Delta q &= q \otimes  q_\Delta  - q \\
	              &=  q \otimes (q_\Delta - q_I)  
\end{align}
Now if we assume that the update is a small rotation so that $q_\Delta = q(\delta a)$ and use the first order approximation from Eq. \eqref{eq:qdela1}
\begin{align}
	\Delta q &= q \otimes  (\delta q (\vect{\delta a}) - q_I) \\
		&= q \otimes  (\begin{bmatrix}1\\ \vect{\delta a}/2\end{bmatrix} - \begin{bmatrix}1\\ \vect{0}\end{bmatrix} ) \\
	\Delta q &= \frac{1}{2} q \otimes \begin{bmatrix}0\\ \vect{\delta a}\end{bmatrix} \label{eq:deltaq}
\end{align}
The similarity of this equation with Eq. \eqref{eq:qdota} demonstrates that if we estimate a small rotation as $\vect{\delta a} = \vect{\omega}\Delta t$, then \eqref{eq:deltaq} is the same as the rectangular rule numerical integration of $\dot{q}$. 

If we now refine further the approximation of the quaternion update using the first two terms of the Taylor's series expansion
\begin{align}
	\Delta q = \dot{q} \Delta t + \frac{1}{2}\ddot{q}\Delta t^2
\end{align}
Then substituting Eq. \eqref{eq:qdota} and its derivative
\begin{align}
	\ddot{q} &= \frac{1}{2} \Biggl( \dot{q} \otimes \begin{bmatrix} 0 \\ \vect{\omega} \end{bmatrix} + q \otimes \begin{bmatrix} 0 \\ \dot{\vect{\omega}} \end{bmatrix} \Biggr)
\end{align}
we reach
\begin{align}
	\Delta q = \frac{1}{2} q \otimes \Biggl( \begin{bmatrix} 0 \\ \vect{\omega} \, \Delta t \end{bmatrix} + 
		\begin{bmatrix} -\Delta t^2 |\vect{\omega}|^2 /4 \\ \dot{\vect{\omega}} \, \Delta t^2 /2\end{bmatrix}   \Biggr)
\end{align}
where the terms are split to show the first order term separated from the second order terms.  The error rate incurred by using only the first order approximation is dominated by $\Delta t |\vect{\omega}|^2 /4$ which may not be negligible for some applications.

%
%
%
%
%



\bibliographystyle{IEEEtran}
\bibliography{qiqqa,books,papers,urls}

%








\end{document}

%% file: estimator_results_nominal/performance_table.tex
\newcolumntype{Y}{>{\raggedleft\arraybackslash}X}
\noindent
\begin{table}
\caption{Results for Scenario "nominal"}
\label{table:nominal}
\begin{center}
\begin{tabularx}{250pt}{| r | c | Y | Y | Y | Y |}
\hline
\textbf{Test case} &  \textbf{Est} & MaxEVz (deg) & MaxEVxy (deg) & FinH (deg) & FinPR (deg) \\ 
 \hline
mockup\_long\_hover &  ekf &    1.24 &    0.58 &    0.35 &    0.38 \\ 
mockup\_long\_hover &   cf &    2.78 &    0.95 &    0.22 &    0.14 \\ 
mockup\_long\_hover &  ukf &    0.89 &    0.78 &    0.20 &    0.52 \\ 
      mockup\_easy &  ekf &    1.80 &    0.65 &    0.02 &    0.28 \\ 
      mockup\_easy &   cf &    3.30 &    1.08 &    0.27 &    0.52 \\ 
      mockup\_easy &  ukf &    0.87 &    0.73 &    0.45 &    0.48 \\ 
   mockup\_slowrot &  ekf &    0.74 &    1.12 &    0.07 &    0.24 \\ 
   mockup\_slowrot &   cf &    1.86 &    1.63 &    0.01 &    0.42 \\ 
   mockup\_slowrot &  ukf &    0.82 &    0.76 &    0.20 &    0.56 \\ 
            mockup &  ekf &    2.93 &    2.17 &    0.34 &    0.37 \\ 
            mockup &   cf &    1.96 &    1.56 &    0.43 &    0.27 \\ 
            mockup &  ukf &    1.08 &    0.76 &    0.03 &    0.51 \\ 
        straightup &  ekf &    1.21 &    0.59 &    0.20 &    0.26 \\ 
        straightup &   cf &    2.50 &    0.86 &    0.03 &    0.57 \\ 
        straightup &  ukf &    1.70 &    1.03 &    0.37 &    0.26 \\ 
      bumpy\_hover &  ekf &    3.52 &    4.23 &    0.23 &    0.27 \\ 
      bumpy\_hover &   cf &   11.81 &   11.67 &    1.82 &    0.74 \\ 
      bumpy\_hover &  ukf &   28.30 &   21.45 &    0.23 &    1.04 \\ 
    straightflight &  ekf &   18.16 &   17.61 &    0.10 &    0.18 \\ 
    straightflight &   cf &    6.87 &    9.28 &    0.78 &    0.73 \\ 
    straightflight &  ukf &   55.54 &   21.26 &    0.19 &    0.06 \\ 
longturn\_bothways &  ekf &   13.86 &   10.45 &    0.31 &    0.27 \\ 
longturn\_bothways &   cf &   11.34 &    8.79 &    1.15 &    0.98 \\ 
longturn\_bothways &  ukf &   24.05 &   27.86 &    0.05 &    2.78 \\ 
        long\_turn &  ekf &   14.36 &   14.06 &    1.24 &    1.15 \\ 
        long\_turn &   cf &   11.10 &   11.29 &    1.48 &    1.87 \\ 
        long\_turn &  ukf &   24.02 &   15.99 &    0.78 &    0.94 \\ 
\hline\end{tabularx}
\end{center}\end{table}

%% file: estimator_results_no_mag/performance_table.tex
\newcolumntype{Y}{>{\raggedleft\arraybackslash}X}
\noindent
\begin{table}
\caption{Results for Scenario "no\_mag"}
\label{table:no_mag}
\begin{center}
\begin{tabularx}{250pt}{| r | c | Y | Y | Y | Y |}
\hline
\textbf{Test case} &  \textbf{Est} & MaxEVz (deg) & MaxEVxy (deg) & FinH (deg) & FinPR (deg) \\ 
 \hline
mockup\_long\_hover &  ekf &    2.23 &    0.58 &    2.23 &    0.38 \\ 
mockup\_long\_hover &   cf &    2.19 &    0.94 &    2.19 &    0.14 \\ 
mockup\_long\_hover &  ukf &    2.21 &    0.40 &    2.17 &    0.17 \\ 
      mockup\_easy &  ekf &    4.92 &    0.64 &    4.92 &    0.30 \\ 
      mockup\_easy &   cf &    4.13 &    0.90 &    4.13 &    0.56 \\ 
      mockup\_easy &  ukf &    4.70 &    0.90 &    4.70 &    0.09 \\ 
   mockup\_slowrot &  ekf &    1.92 &    1.37 &    1.92 &    0.41 \\ 
   mockup\_slowrot &   cf &    3.18 &    1.39 &    3.18 &    0.56 \\ 
   mockup\_slowrot &  ukf &    1.31 &    1.32 &    1.22 &    0.05 \\ 
            mockup &  ekf &    9.20 &    2.12 &    9.20 &    0.43 \\ 
            mockup &   cf &    9.92 &    2.49 &    9.92 &    0.38 \\ 
            mockup &  ukf &    7.43 &    2.13 &    7.43 &    0.09 \\ 
        straightup &  ekf &    4.93 &    0.59 &    4.92 &    0.26 \\ 
        straightup &   cf &    5.08 &    0.86 &    5.07 &    0.58 \\ 
        straightup &  ukf &    6.00 &    1.04 &    5.94 &    0.27 \\ 
      bumpy\_hover &  ekf &   16.08 &   10.87 &   16.08 &    0.32 \\ 
      bumpy\_hover &   cf &    9.91 &    6.95 &    9.91 &    0.61 \\ 
      bumpy\_hover &  ukf &   29.38 &   23.15 &   29.38 &    0.27 \\ 
    straightflight &  ekf &    3.18 &   21.90 &    0.15 &    0.23 \\ 
    straightflight &   cf &    4.14 &    9.14 &    4.14 &    0.19 \\ 
    straightflight &  ukf &   13.03 &   46.18 &   12.70 &    0.33 \\ 
longturn\_bothways &  ekf &    8.49 &   10.30 &    8.33 &    0.62 \\ 
longturn\_bothways &   cf &    6.15 &    6.78 &    6.16 &    1.06 \\ 
longturn\_bothways &  ukf &   31.61 &   18.44 &   31.61 &    0.40 \\ 
        long\_turn &  ekf &   14.73 &   15.27 &   14.03 &    2.03 \\ 
        long\_turn &   cf &    9.07 &   10.11 &    7.92 &    1.77 \\ 
        long\_turn &  ukf &   35.65 &   25.41 &   35.65 &    0.91 \\ 
\hline\end{tabularx}
\end{center}\end{table}

%% file: estimator_results_3D_mag/performance_table.tex
\newcolumntype{Y}{>{\raggedleft\arraybackslash}X}
\noindent
\begin{table}
\caption{Results for Scenario "3D\_mag"}
\label{table:3D_mag}
\begin{center}
\begin{tabularx}{250pt}{| r | c | Y | Y | Y | Y |}
\hline
\textbf{Test case} &  \textbf{Est} & MaxEVz (deg) & MaxEVxy (deg) & FinH (deg) & FinPR (deg) \\ 
 \hline
mockup\_long\_hover &  ekf &    0.98 &    0.73 &    0.04 &    0.19 \\ 
mockup\_long\_hover &   cf &    2.17 &    1.62 &    0.39 &    0.52 \\ 
mockup\_long\_hover &  ukf &    0.52 &    0.42 &    0.09 &    0.26 \\ 
      mockup\_easy &  ekf &    1.15 &    0.62 &    0.03 &    0.17 \\ 
      mockup\_easy &   cf &    2.51 &    1.82 &    0.57 &    0.30 \\ 
      mockup\_easy &  ukf &    0.58 &    0.55 &    0.06 &    0.24 \\ 
   mockup\_slowrot &  ekf &    0.57 &    1.02 &    0.25 &    0.11 \\ 
   mockup\_slowrot &   cf &    1.57 &    1.34 &    0.11 &    0.12 \\ 
   mockup\_slowrot &  ukf &    0.60 &    0.60 &    0.09 &    0.09 \\ 
            mockup &  ekf &    2.96 &    2.07 &    0.89 &    0.64 \\ 
            mockup &   cf &    2.81 &    2.04 &    1.90 &    1.18 \\ 
            mockup &  ukf &    0.65 &    0.60 &    0.16 &    0.18 \\ 
        straightup &  ekf &    0.91 &    0.73 &    0.01 &    0.05 \\ 
        straightup &   cf &    2.38 &    1.80 &    0.98 &    0.52 \\ 
        straightup &  ukf &    0.82 &    0.60 &    0.01 &    0.07 \\ 
      bumpy\_hover &  ekf &    1.87 &    1.34 &    0.46 &    0.14 \\ 
      bumpy\_hover &   cf &    1.77 &    1.90 &    0.78 &    0.29 \\ 
      bumpy\_hover &  ukf &   26.25 &   17.18 &    0.28 &    0.04 \\ 
    straightflight &  ekf &    2.04 &    1.01 &    0.50 &    0.20 \\ 
    straightflight &   cf &    3.03 &    2.89 &    1.00 &    0.37 \\ 
    straightflight &  ukf &   29.86 &   10.59 &    0.16 &    0.01 \\ 
longturn\_bothways &  ekf &    4.26 &    3.52 &    0.66 &    0.39 \\ 
longturn\_bothways &   cf &    5.26 &    4.99 &    0.60 &    0.25 \\ 
longturn\_bothways &  ukf &   17.65 &   16.87 &    0.29 &    0.13 \\ 
        long\_turn &  ekf &    4.50 &    4.20 &    0.96 &    0.50 \\ 
        long\_turn &   cf &    4.38 &    4.23 &    2.29 &    1.24 \\ 
        long\_turn &  ukf &   16.16 &   17.59 &    0.27 &    0.25 \\ 
\hline\end{tabularx}
\end{center}\end{table}

%% file: estimator_results_dynamic_gains/performance_table.tex
\newcolumntype{Y}{>{\raggedleft\arraybackslash}X}
\noindent
\begin{table}
\caption{Results for Scenario "dynamic\_gains"}
\label{table:dynamic_gains}
\begin{center}
\begin{tabularx}{250pt}{| r | c | Y | Y | Y | Y |}
\hline
\textbf{Test case} &  \textbf{Est} & MaxEVz (deg) & MaxEVxy (deg) & FinH (deg) & FinPR (deg) \\ 
 \hline
mockup\_long\_hover &  ekf &    0.71 &    0.41 &    0.08 &    0.14 \\ 
mockup\_long\_hover &   cf &    1.41 &    0.36 &    0.31 &    0.10 \\ 
mockup\_long\_hover &  ukf &    0.93 &    0.44 &    0.14 &    0.19 \\ 
      mockup\_easy &  ekf &    0.59 &    0.35 &    0.07 &    0.10 \\ 
      mockup\_easy &   cf &    1.35 &    0.59 &    0.93 &    0.09 \\ 
      mockup\_easy &  ukf &    0.86 &    0.43 &    0.45 &    0.15 \\ 
   mockup\_slowrot &  ekf &    0.66 &    0.71 &    0.03 &    0.04 \\ 
   mockup\_slowrot &   cf &    1.19 &    1.05 &    0.33 &    0.05 \\ 
   mockup\_slowrot &  ukf &    0.82 &    0.78 &    0.26 &    0.07 \\ 
            mockup &  ekf &    0.69 &    0.52 &    0.03 &    0.04 \\ 
            mockup &   cf &    1.46 &    0.82 &    0.09 &    0.03 \\ 
            mockup &  ukf &    1.02 &    0.65 &    0.03 &    0.10 \\ 
        straightup &  ekf &    1.09 &    0.62 &    0.27 &    0.25 \\ 
        straightup &   cf &    1.65 &    0.44 &    0.41 &    0.22 \\ 
        straightup &  ukf &    1.72 &    1.00 &    0.48 &    0.24 \\ 
      bumpy\_hover &  ekf &    2.90 &    3.75 &    0.12 &    0.15 \\ 
      bumpy\_hover &   cf &   11.99 &   12.43 &    0.45 &    0.10 \\ 
      bumpy\_hover &  ukf &   43.05 &   30.37 &    0.22 &    0.29 \\ 
    straightflight &  ekf &   15.40 &   14.79 &    0.11 &    0.43 \\ 
    straightflight &   cf &    7.12 &    8.19 &    0.45 &    0.45 \\ 
    straightflight &  ukf &   52.02 &   18.37 &    0.28 &    0.26 \\ 
longturn\_bothways &  ekf &   14.97 &   11.50 &    0.29 &    0.10 \\ 
longturn\_bothways &   cf &   11.62 &    8.79 &    0.39 &    0.10 \\ 
longturn\_bothways &  ukf &   23.60 &   27.35 &    0.09 &    0.56 \\ 
        long\_turn &  ekf &   14.08 &   14.06 &    0.40 &    0.49 \\ 
        long\_turn &   cf &   11.10 &   11.29 &    0.12 &    0.40 \\ 
        long\_turn &  ukf &   24.02 &   15.99 &    0.64 &    0.62 \\ 
\hline\end{tabularx}
\end{center}\end{table}